\documentclass[10pt,a4paper]{amsart}

\usepackage[austrian,english]{babel}
\usepackage[latin1]{inputenc}
\usepackage[a4paper]{geometry}
\usepackage{hyperref}
\usepackage{color}
\usepackage{tikz}
\usetikzlibrary{arrows,decorations.pathreplacing,patterns}
\usepackage[font=footnotesize]{caption}

\usepackage{amsthm,amsmath,amssymb}

\newcommand{\arxiv}[1]{\href{http://arxiv.org/abs/#1}{\texttt{arXiv:#1}}}

\allowdisplaybreaks[4]

\newcommand{\mytilde}{\raise.17ex\hbox{$\scriptstyle \mathtt{\sim}$}}

\newcommand{\abs}[1]{\left| #1 \right|}

\newcommand{\ceil}[1]{\left\lceil #1 \right\rceil}

\renewcommand{\d}{\,\mathrm{d}}

\newcommand{\E}{\mathbb{E}}

\newcommand{\N}{\mathbb{N}}
\newcommand{\op}[1]{\operatorname{#1}}

\newcommand{\Q}{\mathbb{Q}}

\newcommand{\R}{\mathbb{R}}

\newcommand{\Z}{\mathbb{Z}}
\renewcommand{\S}{\mathfrak{S}}

\newcommand{\myfontA}[1]{\textnormal{\textbf{#1}}}
\newcommand{\myfontB}[2][]{\myfontA{#2}\textnormal{#1}}

\newcommand{\sk}{

\smallskip}

\newcounter{mythrm}\numberwithin{mythrm}{section}
\numberwithin{equation}{section}

\newenvironment{mymathenvironmentA}[3]{
\refstepcounter{mythrm}\label{#3:#2}
\sk\myfontA{#3 \arabic{section}.\arabic{mythrm}#1.}\begin{itshape}}
{\end{itshape}}

\newenvironment{mymathenvironmentB}[3]{
\refstepcounter{mythrm}\label{#3:#2}
\smallskip\myfontB{#3 \arabic{section}.\arabic{mythrm}#1.}}
{\hfill$\Box$}

\newenvironment{myconj}[2][]{
\begin{mymathenvironmentA}{#1}{#2}{Conjecture}}
{\end{mymathenvironmentA}}

\newenvironment{mylem}[2][]{
\begin{mymathenvironmentA}{#1}{#2}{Lemma}}
{\end{mymathenvironmentA}}

\newenvironment{myproof*}[1]{
\smallskip\myfontA{Proof #1.}}
{\hfill$\blacksquare$}

\newenvironment{myprop}[2][]{
\begin{mymathenvironmentA}{#1}{#2}{Proposition}}
{\end{mymathenvironmentA}}

\newcommand{\refx}[2]{#2~\ref{#2:#1}}

\newcommand{\refcj}[1]{\refx{#1}{Conjecture}}

\newcommand{\reff}[1]{\refx{#1}{Figure}}
\newcommand{\refl}[1]{\refx{#1}{Lemma}}
\newcommand{\refp}[1]{\refx{#1}{Proposition}}
\newcommand{\refq}[1]{\eqref{eq:#1}}

\newcommand{\refs}[1]{\refx{#1}{Section}}
\newcommand{\reft}[1]{\refx{#1}{Theorem}}

\newtheorem{thm}{Theorem}[section]
\newtheorem{conj}{Conjecture}
\newtheorem{lem}[thm]{Lemma}
\newtheorem{prop}[thm]{Proposition}
\newenvironment{mythm}[1]{\begin{thm}\label{Theorem:#1}}{\end{thm}}
\renewenvironment{myconj}[1]{\begin{conj}\label{Conjecture:#1}}{\end{conj}}
\renewenvironment{mylem}[1]{\begin{lem}\label{Lemma:#1}}{\end{lem}}
\renewenvironment{myprop}[1]{\begin{prop}\label{Proposition:#1}}{\end{prop}}

\theoremstyle{remark}
\newtheorem{rem}[thm]{Remark}
\newcommand{\K}{\mathbb{K}}

\newcommand{\arm}{\operatorname{arm}}
\newcommand{\leg}{\operatorname{leg}}

\definecolor{blaugrau}{rgb}{0.796887, 0.789075, 0.871107}

\newcounter{mmacnt}
\def\restartmma{\setcounter{mmacnt}{0}}
\restartmma \catcode`|=\active
\def|#1|{\mathrm{#1}}
\catcode`|=12
\newenvironment{mma}{
 \par
 \catcode`|=\active
 \parskip=2pt\parindent=0pt 
 \small
 \def\In##1\\{%
   \def\linebreak{\hfill\break\null\qquad}%
   \refstepcounter{mmacnt}
   \hangindent=2.5em\hangafter=0
   \leavevmode
   \llap{\tiny\sffamily In[\arabic{mmacnt}]:=\kern.5em}%
   \mathversion{bold}\scriptsize$\tt\bf\displaystyle##1$\normalsize
   \mathversion{normal}\par
 }%
 \def\Print##1\\{%
   \def\linebreak{\hfill\break}%
   \hangindent=2.5em\hangafter=0
   \leavevmode\scriptsize ##1\par}%
 \def\Out##1\\{%
   \vspace*{-0.2cm}\def\linebreak{$\hfill\break\null\hfill$}%
   \kern\abovedisplayskip\par
   \hangindent=2.5em\hangafter=0
   \leavevmode
   \llap{\tiny\sffamily Out[\arabic{mmacnt}]=\kern.5em}
   \scriptsize$\displaystyle\tt##1$\normalsize\hfill\null\par
   \kern\belowdisplayskip\vspace*{-0.1cm}
 }%
 \def\Warning##1##2\\{%
   \def\linebreak{\hfill\break}%
   \hangindent=2.5em\hangafter=0
   \leavevmode
   {\scriptsize##1 : ##2}\par}%
}{%
 \par\smallskip
}

\newcommand{\LoadP}[1]{\fcolorbox{black}{blaugrau}{
\begin{minipage}[t]{13.4cm}
\footnotesize #1
\end{minipage}}}

\newcommand{\myOut}[1]{{\sffamily Out[#1]}}

\def\MLabel#1{{\refstepcounter{mmacnt}\label{#1}}\addtocounter{mmacnt}{-1}}

\begin{document}


\title[Asymptotic and exact results on the complexity of the NPS algorithm]{Asymptotic and exact results on the complexity of the Novelli--Pak--Stoyanovskii algorithm}

\author{Carsten Schneider}
\address[Carsten Schneider]{Research Institute for Symbolic Computation\\
J. Kepler University Linz\\
A-4040 Linz, Austria}
\email{Carsten.Schneider@risc.jku.at}

\author{Robin Sulzgruber}
\address[Robin Sulzgruber]{Fakult\"at f\"ur Mathematik, Universit\"at Wien, Oskar-Morgenstern-Platz 1, 1090 Wien, Austria}
\email{robin.sulzgruber@univie.ac.at}

\thanks{Research supported by the Austrian Science Fund (FWF), grants F5005-N15 and F5009-N15 in the framework of the Special Research Program ``Algorithmic and Enumerative Combinatorics'' (SFB F50).}
\date{May, 2017}

\begin{abstract}
\noindent The Novelli--Pak--Stoyanovskii algorithm is a sorting algorithm for Young  of a fixed shape that was originally devised to give a bijective proof of the hook-length formula.
We obtain new asymptotic results on the average case and worst case complexity of this algorithm as the underlying shape tends to a fixed limit curve.
Furthermore, using the summation package Sigma
we prove an exact formula for the average case complexity when the underlying shape consists of only two rows. 
We thereby answer questions posed by Krattenthaler and M{\"u}ller.

\end{abstract}

\thispagestyle{empty}
\maketitle




\pagenumbering{arabic}
\pagestyle{headings}

The Novelli--Pak--Stoyanovskii algorithm (NPS algorithm) transforms (sorts) an arbitrary filling of a Young diagram $\lambda$ of size $n$ with the numbers $\{1,2,\dots,n\}$ into a standard Young tableau of the shape $\lambda$.
From a sorters point of view, the algorithm is best described as a ``two-dimensional insertion sort''.
Following certain rules, at each step two entries of adjacent cells are compared and possibly exchanged.
The remarkable property of the NPS algorithm is that, when applied to all possible fillings of a fixed diagram, it produces every standard Young tableau of this shape equally often as an output.
Thus it provides a uniformly distributed random sampler for standard Young tableaux of a given shape.
The algorithm was originally defined by Novelli, Pak and Stoyanovskii in \cite{NPS1997,PakSto1992} to give a bijective proof of the hook-length formula due to Frame, Robinson and Thrall \cite{FRT1954}.

The study of the complexity of the NPS algorithm on a partition $\lambda$ was initiated by Krattenthaler and M{\"u}ller.
They define the average case complexity $C(\lambda)$ and the worst case complexity $W(\lambda)$ as the average, respectively the maximal number of exchanges performed by the NPS algorithm applied to an arbitrary filling of shape $\lambda$.
Their analysis lead to multiple interesting conjectures that are the main motivation for the present paper.
Three of these conjectures, which were presented to us in private communication~\cite{Kratt2013}, are listed below.
The first two conjectures concern the asymptotic behavior of $C(\lambda^{(n)})$ and $W(\lambda^{(n)})$, where $(\lambda^{(n)})_{n\in\N}$ is a sequence of partitions that 
approach a fixed limit curve $\gamma$ after a suitable rescaling.
(This is made precise in \refs{hookcoord}.)

\begin{myconj}{C} The order of magnitude of $C(\lambda^{(n)})$ lies between $\abs{\lambda^{(n)}}^{3/2}$ and $\abs{\lambda^{(n)}}^2$ where $\abs{\lambda^{(n)}}$ denotes the size of $\lambda^{(n)}$.
\end{myconj}

\begin{myconj}{CW} The average case complexity $C(\lambda^{(n)})$ is asymptotically one half of the worst case complexity $W(\lambda^{(n)})$.
\end{myconj}

The third conjecture that we are interested in is an exact formula for the average case complexity when the partition $\lambda=(\lambda_1,\lambda_2)$ consists of two parts.

\begin{myconj}{zweizeiler} The average case complexity of the NPS algorithm on $\lambda=(\lambda_1,\lambda_2)$ is given by\footnote{Here $(x)_k$ stands for the Pochhammer symbol which is $1$ for $k=0$ and which equals $x(x+1)\dots(x+k-1)$ for positive integers $k$.
}
\begin{align*}
C(\lambda)
&=\frac{\lambda_1(\lambda_1-1)}{4}+\frac{\lambda_2(\lambda_2-3)}{4}-2\sum_{k=1}^{\lambda_2}\binom{\lambda_2}{k}\frac{(-1)^k(2k-2)!}{(\lambda_1-\lambda_2+2)_{2k-1}}.
\end{align*}
\end{myconj}

In this paper we prove Conjectures~\ref{Conjecture:C} and~\ref{Conjecture:zweizeiler}.
While we were unable to prove \refcj{CW} in full generality, we do provide a proof for a large class of sequences of partitions, adding further evidence to its validity.
The article is structured as follows:

In \refs{algorithm} we review the combinatorics of the NPS algorithm.

In Sections~\ref{Section:hookcoord}--\ref{Section:C} we engage sequences of partitions that converge under a balanced scaling, that is, by a factor of $\sqrt{n}$ in both dimensions.
\refs{hookcoord} contains mostly preparatory results and includes a precise definition of convergence.
\refs{W} treats the worst case complexity.
We derive an exact combinatorial formula for $W(\lambda)$ for any partition $\lambda$ in \refp{W} by proving that a trivial upper bound is tight.
Moreover \reft{Wn} provides an asymptotic result on $W(\lambda^{(n)})$ in the balanced case.
\refs{C} treats the average case complexity in the balanced case.
Here we give an asymptotic lower bound for $C(\lambda^{(n)})$ in \reft{Cn}.
It is a consequence of our results that both $C(\lambda^{(n)})$ and $W(\lambda^{(n)})$ are of order $\abs{\lambda^{(n)}}^{3/2}$, which is in accordance with \refcj{C}.

In \refs{pq} we turn to sequences that converge when subjected to an imbalanced scaling, that is, by a factor of $n^{1/p}$ in one direction and a factor of $n^{1/q}$ in another.
\reft{CWpq} verifies both Conjectures~\ref{Conjecture:C} and~\ref{Conjecture:CW} in the imbalanced setting.
More precisely, we show that $C(\lambda^{(n)})$ and $W(\lambda^{(n)})$ are both of order $n^{(p+1)/p}$ if $p<q$ and of order $n^{(q+1)/q}$ if $p>q$, and that the leading coefficient of the average case complexity is one half of the leading coefficient of the worst case complexity.

In \refs{zweizeiler} we prove \refcj{zweizeiler} employing the summation package \texttt{Sigma}~\cite{Sigma} in a non-trivial fashion. Here we first provide an alternative representation $C(\lambda)$ in terms of five non-trivial double sums and show that this expression corresponds to the single sum expression given in \refcj{zweizeiler}. The underlying machinery is based on the summation paradigms of creative telescoping, recurrence solving and the zero-recognition problem for the class of (indefinite) nested sums over hypergeometric products. As a by-product we provide alternative representations of $C(\lambda)=C(\lambda_1,\lambda_2)$ that enable one to calculate $C(\lambda)$ efficiently if one keeps $\lambda_1$ symbolic and specializes $\lambda_2$ to a concrete value, or if one keeps $\lambda_1$ symbolic and specializes the distance $\lambda_1-\lambda_2$ to a specific non-negative integer. In particular, we discover a particularly nice formula for the special case $\lambda_1=\lambda_2$.



\section{The NPS algorithm}\label{Section:algorithm}

\begin{figure}[ht]
\begin{center}
\begin{tikzpicture}[scale=.5]
\draw
(0,1)--(1,1)
(0,2)--(1,2)
(0,3)--(1,3)
(0,4)--(2,4)
(0,5)--(4,5)
(1,3)--(1,6)
(2,4)--(2,6)
(3,4)--(3,6);
\draw[thick] (0,0)--(0,6)--(4,6)--(4,4)--(2,4)--(2,3)--(1,3)--(1,0)--cycle;
\begin{scope}[xshift=10cm]
\draw
(0,1)--(1,1)
(0,2)--(1,2)
(0,3)--(1,3)
(0,4)--(2,4)
(0,5)--(4,5)
(1,3)--(1,6)
(2,4)--(2,6)
(3,4)--(3,6);
\draw[thick] (0,0)--(0,6)--(4,6)--(4,4)--(2,4)--(2,3)--(1,3)--(1,0)--cycle;
\draw[xshift=5mm,yshift=5mm]
(0,5)node{$12$}
(1,5)node{$7$}
(2,5)node{$5$}
(3,5)node{$1$}
(0,4)node{$2$}
(1,4)node{$10$}
(2,4)node{$9$}
(3,4)node{$11$}
(0,3)node{$13$}
(1,3)node{$4$}
(0,2)node{$8$}
(0,1)node{$6$}
(0,0)node{$3$};
\end{scope}
\end{tikzpicture}
\caption{The partition $\lambda=(4,4,2,1,1,1)$ in English convention, and a tableau of shape $\lambda$.}
\label{Figure:english}
\end{center}
\end{figure}
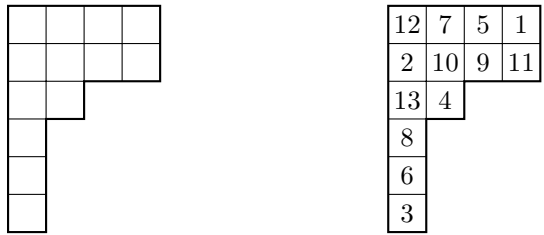

In this section we recall some definitions concerning partitions and Young tableaux as well as the needed facts about the NPS algorithm.

\smallskip
Let $n\in\N$ be a non-negative integer.
A \emph{partition} $\lambda$ of $n$ is a weakly decreasing sequence $\lambda_1\geq\lambda_2\geq\dots\geq\lambda_k>0$ of positive integers such that $\abs{\lambda}=\sum\lambda_i=n$.
We call $\abs{\lambda}$ the \emph{size} of $\lambda$.
The \emph{length} $l(\lambda)$ is the number of summands $\lambda_i$.
We identify a partition with its \emph{Young diagram} $\lambda=\{(i,j):1\leq i\leq l(\lambda),1\leq j\leq \lambda_i\}$.
The elements $(i,j)$ are called the \emph{cells} of the partition $\lambda$.

The \emph{conjugate} partition $\lambda'$ of $\lambda$ is given by the Young diagram $\{(j,i):(i,j)\in\lambda\}$.
Visually we imagine a partition as a left justified array of $n$ boxes, with $\lambda_i$ boxes in the $i$-th row counting from top to bottom as in \reff{english}.
Thus $\lambda'$ is obtained from $\lambda$ by a flip along the main diagonal.


Define the \emph{arm} of a cell $\arm_{\lambda}(i,j)=\lambda_i-j$ as the number of cells in the same row as $(i,j)$ and strictly to the right of $(i,j)$.
Define the \emph{leg} of a cell $\leg_{\lambda}(i,j)=\lambda_j'-i$ as the number of cells in the same column as $(i,j)$ and strictly below $(i,j)$.
Furthermore, define the \emph{hook length} of a cell as $h_{\lambda}(i,j)=\lambda_i+\lambda_j'-i-j+1$.
We call a cell $(i,j)$ a \emph{corner} of $\lambda$ if $h_{\lambda}(i,j)=1$.

An \emph{integer filling} of a partition is a map $T:\lambda\to\Z$ assigning an entry $T(i,j)$ to each cell $(i,j)$.
The partition $\lambda$ is called the \emph{shape} of $T$.
A \emph{tableau} is a bijection $T:\lambda\to\{1,2,\dots,n\}$.
A \emph{standard Young tableau} (SYT) is a tableau that increases along rows from left to right as well as down columns, that is, $T(i,j)\leq T(i',j')$ whenever $i\leq i'$ and $j\leq j'$.
A \emph{hook tableau} is a map $H:\lambda\to\Z$ such that for each cell $(i,j)$ we have $-\leg_{\lambda}(i,j)\leq H(i,j)\leq\arm_{\lambda}(i,j)$.
Note that somewhat counter-intuitively a hook tableau is not a tableau.
We denote the set of all tableaux, standard Young tableaux and hook tableaux of shape $\lambda$ by $\op{T}(\lambda)$, $\op{SYT}(\lambda)$ and $\op{H}(\lambda)$, respectively.

\smallskip
The celebrated hook-length formula~\cite{FRT1954} determines the number of standard Young tableaux of a fixed shape:
\begin{align}\label{eq:hlf}
\#\op{SYT}(\lambda) 
=\frac{n!}{\displaystyle\prod_{(i,j)\in\lambda}h_{\lambda}(i,j)}.
\end{align}
Since $\#\op{T}(\lambda)=n!$ and $\#\op{H}(\lambda)=\prod_{(i,j)\in\lambda}h_{\lambda}(i,j)$ it is possible to prove \refq{hlf} bijectively by constructing a bijection
\begin{align*}
\Phi:\op{T}(\lambda) \to \op{SYT}(\lambda) \times \op{H}(\lambda).
\end{align*}
Such a construction was found by Novelli, Pak and Stoyanovskii \cite{PakSto1992,NPS1997}.
We are now going to describe the map $\Phi$.
See \reff{nps} for an example.

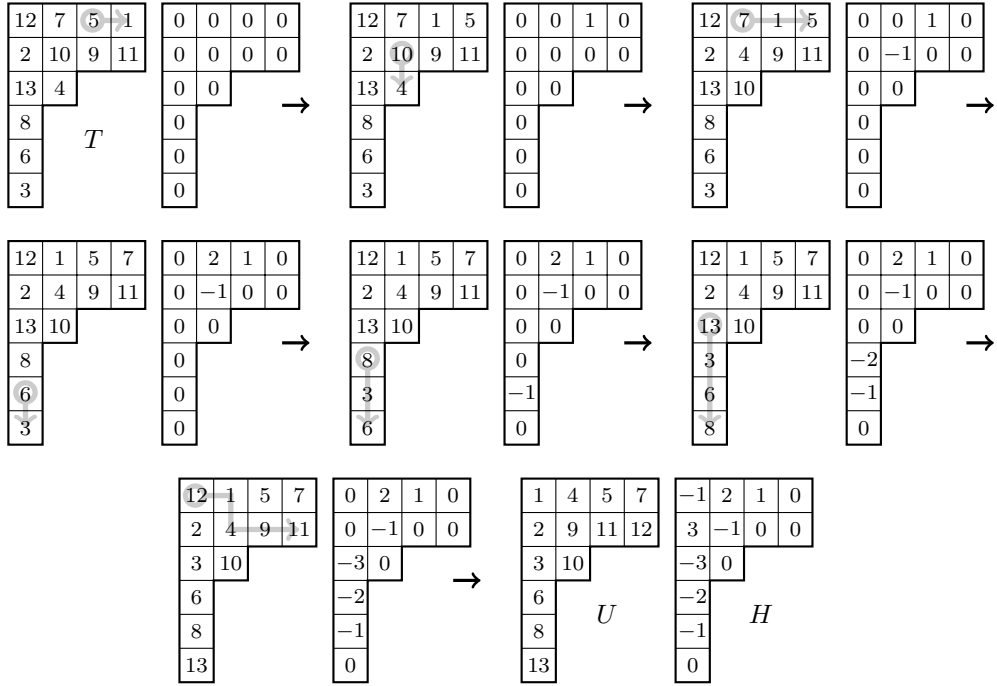
\begin{figure}[t]
\begin{center}
\begin{tikzpicture}[scale=.45]
\draw[black!20,line width=2pt,o->,xshift=5mm,yshift=5mm](1.6,5)--(3,5);
\begin{scope}[xshift=0cm]
\draw[very thick,->](8,3)--(8.8,3);
\draw(2.5,2)node{$T$};
\draw (0,1)--(1,1) (0,2)--(1,2) (0,3)--(1,3) (0,4)--(2,4) (0,5)--(4,5) (1,3)--(1,6) (2,4)--(2,6) (3,4)--(3,6);
\draw[thick] (0,0)--(0,6)--(4,6)--(4,4)--(2,4)--(2,3)--(1,3)--(1,0)--cycle;
\draw[xshift=5mm,yshift=5mm]
(0,5)node{\footnotesize{$12$}} (1,5)node{\footnotesize{$7$}} (2,5)node{\footnotesize{$5$}} (3,5)node{\footnotesize{$1$}}
(0,4)node{\footnotesize{$2$}} (1,4)node{\footnotesize{$10$}} (2,4)node{\footnotesize{$9$}} (3,4)node{\footnotesize{$11$}}
(0,3)node{\footnotesize{$13$}} (1,3)node{\footnotesize{$4$}}
(0,2)node{\footnotesize{$8$}}
(0,1)node{\footnotesize{$6$}}
(0,0)node{\footnotesize{$3$}};
\begin{scope}[xshift=45mm]
\draw (0,1)--(1,1) (0,2)--(1,2) (0,3)--(1,3) (0,4)--(2,4) (0,5)--(4,5) (1,3)--(1,6) (2,4)--(2,6) (3,4)--(3,6);
\draw[thick] (0,0)--(0,6)--(4,6)--(4,4)--(2,4)--(2,3)--(1,3)--(1,0)--cycle;
\draw[xshift=5mm,yshift=5mm]
(0,5)node{\footnotesize{$0$}} (1,5)node{\footnotesize{$0$}} (2,5)node{\footnotesize{$0$}} (3,5)node{\footnotesize{$0$}}
(0,4)node{\footnotesize{$0$}} (1,4)node{\footnotesize{$0$}} (2,4)node{\footnotesize{$0$}} (3,4)node{\footnotesize{$0$}}
(0,3)node{\footnotesize{$0$}} (1,3)node{\footnotesize{$0$}}
(0,2)node{\footnotesize{$0$}}
(0,1)node{\footnotesize{$0$}}
(0,0)node{\footnotesize{$0$}};
\end{scope}
\end{scope}
\begin{scope}[xshift=10cm]
\draw[black!20,line width=2pt,o->,xshift=5mm,yshift=5mm](1,4.4)--(1,3);
\draw[very thick,->](8,3)--(8.8,3);
\draw (0,1)--(1,1) (0,2)--(1,2) (0,3)--(1,3) (0,4)--(2,4) (0,5)--(4,5) (1,3)--(1,6) (2,4)--(2,6) (3,4)--(3,6);
\draw[thick] (0,0)--(0,6)--(4,6)--(4,4)--(2,4)--(2,3)--(1,3)--(1,0)--cycle;
\draw[xshift=5mm,yshift=5mm]
(0,5)node{\footnotesize{$12$}} (1,5)node{\footnotesize{$7$}} (2,5)node{\footnotesize{$1$}} (3,5)node{\footnotesize{$5$}}
(0,4)node{\footnotesize{$2$}} (1,4)node{\footnotesize{$10$}} (2,4)node{\footnotesize{$9$}} (3,4)node{\footnotesize{$11$}}
(0,3)node{\footnotesize{$13$}} (1,3)node{\footnotesize{$4$}}
(0,2)node{\footnotesize{$8$}}
(0,1)node{\footnotesize{$6$}}
(0,0)node{\footnotesize{$3$}};
\begin{scope}[xshift=45mm]
\draw (0,1)--(1,1) (0,2)--(1,2) (0,3)--(1,3) (0,4)--(2,4) (0,5)--(4,5) (1,3)--(1,6) (2,4)--(2,6) (3,4)--(3,6);
\draw[thick] (0,0)--(0,6)--(4,6)--(4,4)--(2,4)--(2,3)--(1,3)--(1,0)--cycle;
\draw[xshift=5mm,yshift=5mm]
(0,5)node{\footnotesize{$0$}} (1,5)node{\footnotesize{$0$}} (2,5)node{\footnotesize{$1$}} (3,5)node{\footnotesize{$0$}}
(0,4)node{\footnotesize{$0$}} (1,4)node{\footnotesize{$0$}} (2,4)node{\footnotesize{$0$}} (3,4)node{\footnotesize{$0$}}
(0,3)node{\footnotesize{$0$}} (1,3)node{\footnotesize{$0$}}
(0,2)node{\footnotesize{$0$}}
(0,1)node{\footnotesize{$0$}}
(0,0)node{\footnotesize{$0$}};
\end{scope}
\end{scope}
\begin{scope}[xshift=20cm]
\draw[black!20,line width=2pt,o->,xshift=5mm,yshift=5mm](.6,5)--(3,5);
\draw[very thick,->](8,3)--(8.8,3);
\draw (0,1)--(1,1) (0,2)--(1,2) (0,3)--(1,3) (0,4)--(2,4) (0,5)--(4,5) (1,3)--(1,6) (2,4)--(2,6) (3,4)--(3,6);
\draw[thick] (0,0)--(0,6)--(4,6)--(4,4)--(2,4)--(2,3)--(1,3)--(1,0)--cycle;
\draw[xshift=5mm,yshift=5mm]
(0,5)node{\footnotesize{$12$}} (1,5)node{\footnotesize{$7$}} (2,5)node{\footnotesize{$1$}} (3,5)node{\footnotesize{$5$}}
(0,4)node{\footnotesize{$2$}} (1,4)node{\footnotesize{$4$}} (2,4)node{\footnotesize{$9$}} (3,4)node{\footnotesize{$11$}}
(0,3)node{\footnotesize{$13$}} (1,3)node{\footnotesize{$10$}}
(0,2)node{\footnotesize{$8$}}
(0,1)node{\footnotesize{$6$}}
(0,0)node{\footnotesize{$3$}};
\begin{scope}[xshift=45mm]
\draw (0,1)--(1,1) (0,2)--(1,2) (0,3)--(1,3) (0,4)--(2,4) (0,5)--(4,5) (1,3)--(1,6) (2,4)--(2,6) (3,4)--(3,6);
\draw[thick] (0,0)--(0,6)--(4,6)--(4,4)--(2,4)--(2,3)--(1,3)--(1,0)--cycle;
\draw[xshift=5mm,yshift=5mm]
(0,5)node{\footnotesize{$0$}} (1,5)node{\footnotesize{$0$}} (2,5)node{\footnotesize{$1$}} (3,5)node{\footnotesize{$0$}}
(0,4)node{\footnotesize{$0$}} (1,4)node{\footnotesize{$-1$}} (2,4)node{\footnotesize{$0$}} (3,4)node{\footnotesize{$0$}}
(0,3)node{\footnotesize{$0$}} (1,3)node{\footnotesize{$0$}}
(0,2)node{\footnotesize{$0$}}
(0,1)node{\footnotesize{$0$}}
(0,0)node{\footnotesize{$0$}};
\end{scope}
\end{scope}
\begin{scope}[xshift=0cm,yshift=-7cm]
\draw[black!20,line width=2pt,o->,xshift=5mm,yshift=5mm](0,1.4)--(0,0);
\draw[very thick,->](8,3)--(8.8,3);
\draw (0,1)--(1,1) (0,2)--(1,2) (0,3)--(1,3) (0,4)--(2,4) (0,5)--(4,5) (1,3)--(1,6) (2,4)--(2,6) (3,4)--(3,6);
\draw[thick] (0,0)--(0,6)--(4,6)--(4,4)--(2,4)--(2,3)--(1,3)--(1,0)--cycle;
\draw[xshift=5mm,yshift=5mm]
(0,5)node{\footnotesize{$12$}} (1,5)node{\footnotesize{$1$}} (2,5)node{\footnotesize{$5$}} (3,5)node{\footnotesize{$7$}}
(0,4)node{\footnotesize{$2$}} (1,4)node{\footnotesize{$4$}} (2,4)node{\footnotesize{$9$}} (3,4)node{\footnotesize{$11$}}
(0,3)node{\footnotesize{$13$}} (1,3)node{\footnotesize{$10$}}
(0,2)node{\footnotesize{$8$}}
(0,1)node{\footnotesize{$6$}}
(0,0)node{\footnotesize{$3$}};
\begin{scope}[xshift=45mm]
\draw (0,1)--(1,1) (0,2)--(1,2) (0,3)--(1,3) (0,4)--(2,4) (0,5)--(4,5) (1,3)--(1,6) (2,4)--(2,6) (3,4)--(3,6);
\draw[thick] (0,0)--(0,6)--(4,6)--(4,4)--(2,4)--(2,3)--(1,3)--(1,0)--cycle;
\draw[xshift=5mm,yshift=5mm]
(0,5)node{\footnotesize{$0$}} (1,5)node{\footnotesize{$2$}} (2,5)node{\footnotesize{$1$}} (3,5)node{\footnotesize{$0$}}
(0,4)node{\footnotesize{$0$}} (1,4)node{\footnotesize{$-1$}} (2,4)node{\footnotesize{$0$}} (3,4)node{\footnotesize{$0$}}
(0,3)node{\footnotesize{$0$}} (1,3)node{\footnotesize{$0$}}
(0,2)node{\footnotesize{$0$}}
(0,1)node{\footnotesize{$0$}}
(0,0)node{\footnotesize{$0$}};
\end{scope}
\end{scope}
\begin{scope}[xshift=10cm,yshift=-7cm]
\draw[black!20,line width=2pt,o->,xshift=5mm,yshift=5mm](0,2.4)--(0,0);
\draw[very thick,->](8,3)--(8.8,3);
\draw (0,1)--(1,1) (0,2)--(1,2) (0,3)--(1,3) (0,4)--(2,4) (0,5)--(4,5) (1,3)--(1,6) (2,4)--(2,6) (3,4)--(3,6);
\draw[thick] (0,0)--(0,6)--(4,6)--(4,4)--(2,4)--(2,3)--(1,3)--(1,0)--cycle;
\draw[xshift=5mm,yshift=5mm]
(0,5)node{\footnotesize{$12$}} (1,5)node{\footnotesize{$1$}} (2,5)node{\footnotesize{$5$}} (3,5)node{\footnotesize{$7$}}
(0,4)node{\footnotesize{$2$}} (1,4)node{\footnotesize{$4$}} (2,4)node{\footnotesize{$9$}} (3,4)node{\footnotesize{$11$}}
(0,3)node{\footnotesize{$13$}} (1,3)node{\footnotesize{$10$}}
(0,2)node{\footnotesize{$8$}}
(0,1)node{\footnotesize{$3$}}
(0,0)node{\footnotesize{$6$}};
\begin{scope}[xshift=45mm]
\draw (0,1)--(1,1) (0,2)--(1,2) (0,3)--(1,3) (0,4)--(2,4) (0,5)--(4,5) (1,3)--(1,6) (2,4)--(2,6) (3,4)--(3,6);
\draw[thick] (0,0)--(0,6)--(4,6)--(4,4)--(2,4)--(2,3)--(1,3)--(1,0)--cycle;
\draw[xshift=5mm,yshift=5mm]
(0,5)node{\footnotesize{$0$}} (1,5)node{\footnotesize{$2$}} (2,5)node{\footnotesize{$1$}} (3,5)node{\footnotesize{$0$}}
(0,4)node{\footnotesize{$0$}} (1,4)node{\footnotesize{$-1$}} (2,4)node{\footnotesize{$0$}} (3,4)node{\footnotesize{$0$}}
(0,3)node{\footnotesize{$0$}} (1,3)node{\footnotesize{$0$}}
(0,2)node{\footnotesize{$0$}}
(0,1)node{\footnotesize{$-1$}}
(0,0)node{\footnotesize{$0$}};
\end{scope}
\end{scope}
\begin{scope}[xshift=20cm,yshift=-7cm]
\draw[black!20,line width=2pt,o->,xshift=5mm,yshift=5mm](0,3.4)--(0,0);
\draw[very thick,->](8,3)--(8.8,3);
\draw (0,1)--(1,1) (0,2)--(1,2) (0,3)--(1,3) (0,4)--(2,4) (0,5)--(4,5) (1,3)--(1,6) (2,4)--(2,6) (3,4)--(3,6);
\draw[thick] (0,0)--(0,6)--(4,6)--(4,4)--(2,4)--(2,3)--(1,3)--(1,0)--cycle;
\draw[xshift=5mm,yshift=5mm]
(0,5)node{\footnotesize{$12$}} (1,5)node{\footnotesize{$1$}} (2,5)node{\footnotesize{$5$}} (3,5)node{\footnotesize{$7$}}
(0,4)node{\footnotesize{$2$}} (1,4)node{\footnotesize{$4$}} (2,4)node{\footnotesize{$9$}} (3,4)node{\footnotesize{$11$}}
(0,3)node{\footnotesize{$13$}} (1,3)node{\footnotesize{$10$}}
(0,2)node{\footnotesize{$3$}}
(0,1)node{\footnotesize{$6$}}
(0,0)node{\footnotesize{$8$}};
\begin{scope}[xshift=45mm]
\draw (0,1)--(1,1) (0,2)--(1,2) (0,3)--(1,3) (0,4)--(2,4) (0,5)--(4,5) (1,3)--(1,6) (2,4)--(2,6) (3,4)--(3,6);
\draw[thick] (0,0)--(0,6)--(4,6)--(4,4)--(2,4)--(2,3)--(1,3)--(1,0)--cycle;
\draw[xshift=5mm,yshift=5mm]
(0,5)node{\footnotesize{$0$}} (1,5)node{\footnotesize{$2$}} (2,5)node{\footnotesize{$1$}} (3,5)node{\footnotesize{$0$}}
(0,4)node{\footnotesize{$0$}} (1,4)node{\footnotesize{$-1$}} (2,4)node{\footnotesize{$0$}} (3,4)node{\footnotesize{$0$}}
(0,3)node{\footnotesize{$0$}} (1,3)node{\footnotesize{$0$}}
(0,2)node{\footnotesize{$-2$}}
(0,1)node{\footnotesize{$-1$}}
(0,0)node{\footnotesize{$0$}};
\end{scope}
\end{scope}
\begin{scope}[xshift=5cm,yshift=-14cm]
\draw[black!20,line width=2pt,o->,xshift=5mm,yshift=5mm](-.4,5)--(1,5)--(1,4)--(3,4);
\draw[very thick,->](8,3)--(8.8,3);
\draw (0,1)--(1,1) (0,2)--(1,2) (0,3)--(1,3) (0,4)--(2,4) (0,5)--(4,5) (1,3)--(1,6) (2,4)--(2,6) (3,4)--(3,6);
\draw[thick] (0,0)--(0,6)--(4,6)--(4,4)--(2,4)--(2,3)--(1,3)--(1,0)--cycle;
\draw[xshift=5mm,yshift=5mm]
(0,5)node{\footnotesize{$12$}} (1,5)node{\footnotesize{$1$}} (2,5)node{\footnotesize{$5$}} (3,5)node{\footnotesize{$7$}}
(0,4)node{\footnotesize{$2$}} (1,4)node{\footnotesize{$4$}} (2,4)node{\footnotesize{$9$}} (3,4)node{\footnotesize{$11$}}
(0,3)node{\footnotesize{$3$}} (1,3)node{\footnotesize{$10$}}
(0,2)node{\footnotesize{$6$}}
(0,1)node{\footnotesize{$8$}}
(0,0)node{\footnotesize{$13$}};
\begin{scope}[xshift=45mm]
\draw (0,1)--(1,1) (0,2)--(1,2) (0,3)--(1,3) (0,4)--(2,4) (0,5)--(4,5) (1,3)--(1,6) (2,4)--(2,6) (3,4)--(3,6);
\draw[thick] (0,0)--(0,6)--(4,6)--(4,4)--(2,4)--(2,3)--(1,3)--(1,0)--cycle;
\draw[xshift=5mm,yshift=5mm]
(0,5)node{\footnotesize{$0$}} (1,5)node{\footnotesize{$2$}} (2,5)node{\footnotesize{$1$}} (3,5)node{\footnotesize{$0$}}
(0,4)node{\footnotesize{$0$}} (1,4)node{\footnotesize{$-1$}} (2,4)node{\footnotesize{$0$}} (3,4)node{\footnotesize{$0$}}
(0,3)node{\footnotesize{$-3$}} (1,3)node{\footnotesize{$0$}}
(0,2)node{\footnotesize{$-2$}}
(0,1)node{\footnotesize{$-1$}}
(0,0)node{\footnotesize{$0$}};
\end{scope}
\end{scope}
\begin{scope}[xshift=15cm,yshift=-14cm]
\draw(2.5,2)node{$U$};
\draw (0,1)--(1,1) (0,2)--(1,2) (0,3)--(1,3) (0,4)--(2,4) (0,5)--(4,5) (1,3)--(1,6) (2,4)--(2,6) (3,4)--(3,6);
\draw[thick] (0,0)--(0,6)--(4,6)--(4,4)--(2,4)--(2,3)--(1,3)--(1,0)--cycle;
\draw[xshift=5mm,yshift=5mm]
(0,5)node{\footnotesize{$1$}} (1,5)node{\footnotesize{$4$}} (2,5)node{\footnotesize{$5$}} (3,5)node{\footnotesize{$7$}}
(0,4)node{\footnotesize{$2$}} (1,4)node{\footnotesize{$9$}} (2,4)node{\footnotesize{$11$}} (3,4)node{\footnotesize{$12$}}
(0,3)node{\footnotesize{$3$}} (1,3)node{\footnotesize{$10$}}
(0,2)node{\footnotesize{$6$}}
(0,1)node{\footnotesize{$8$}}
(0,0)node{\footnotesize{$13$}};
\begin{scope}[xshift=45mm]
\draw(2.5,2)node{$H$};
\draw (0,1)--(1,1) (0,2)--(1,2) (0,3)--(1,3) (0,4)--(2,4) (0,5)--(4,5) (1,3)--(1,6) (2,4)--(2,6) (3,4)--(3,6);
\draw[thick] (0,0)--(0,6)--(4,6)--(4,4)--(2,4)--(2,3)--(1,3)--(1,0)--cycle;
\draw[xshift=5mm,yshift=5mm]
(0,5)node{\footnotesize{$-1$}} (1,5)node{\footnotesize{$2$}} (2,5)node{\footnotesize{$1$}} (3,5)node{\footnotesize{$0$}}
(0,4)node{\footnotesize{$3$}} (1,4)node{\footnotesize{$-1$}} (2,4)node{\footnotesize{$0$}} (3,4)node{\footnotesize{$0$}}
(0,3)node{\footnotesize{$-3$}} (1,3)node{\footnotesize{$0$}}
(0,2)node{\footnotesize{$-2$}}
(0,1)node{\footnotesize{$-1$}}
(0,0)node{\footnotesize{$0$}};
\end{scope}
\end{scope}
\end{tikzpicture}
\caption{The NPS algorithm in action.}
\label{Figure:nps}
\end{center}
\end{figure}

First note that given a permutation $\sigma\in\S_n$ and a tableau $T\in\operatorname{T}(\lambda)$, we obtain a new tableau $T'=\sigma\circ T$ by setting $T'(i,j)=\sigma(T(i,j))$.
In particular if $\sigma=(k,m)$ is a transposition, then $(k,m)\circ T$ arises from $T$ by exchanging the two entries $k$ and $m$.

Impose the reverse lexicographic order $\preceq$ on the cells of $\lambda$, that is, $(i,j)\preceq(i',j')$ if $j<j'$ or if $j=j'$ and $i\leq i'$.
Let $(i_1,j_1)\succ(i_2,j_2)\succ\dots\succ(i_n,j_n)$ be the cells of $\lambda$ in decreasing order and set $k_r=T(i_r,j_r)$.
Moreover, set $T_0=T$ and let $H_0$ be the hook tableau with all entries equal to zero.

Given a pair tableau $(T_{r-1},H_{r-1})$ we first construct the tableau $T_r$ from $T_{r-1}$ as follows: 
\begin{description}
\item[E0] Set $T=T_{r-1}$.
\item[E1] Set $(i,j)=T^{-1}(k_r)$.
\item[E2] If $(i,j)$ is a corner of $\lambda$, then return $T_r=T$.
\item[E3] Otherwise set $m=\min T^{-1}(\lambda\cap\{(i+1,j),(i,j+1)\})$.
\item[E4] If $k_r<m$ then return $T_r=T$.
\item[E5] Otherwise $m<k_r$. In this case exchange the entries $m$ and $k_r$, that is, set $T=(k_r,m)\circ T$, and return to step \textbf{E1}.
\end{description}
Note that $T_r$ is obtained from $T_{r-1}$ by applying a jeu-de-taquin-like move to the entry $k_r$.
Next we construct $H_r$ from $H_{r-1}$ as follows:
\begin{description}
\item[H0] Set $H=H_{r-1}$, $(i,j)=T_{r-1}^{-1}(k_r)$ and $(i',j')=T_r^{-1}(k_r)$.
\item[H1\label{item:H1}] For each $s=i,\dots,i'-1$ set $H(s,j)=H_{r-1}(s+1,j)-1$.
\item[H2\label{item:H2}] Set $H(i',j)=j'-j$.
\item[H3] Return $H_r=H$.
\end{description}
These rules give rise to a sequence $(T_0,H_0),(T_1,H_1),\dots,(T_n,H_n)$ of pairs of a tableau and a hook tableau.
We define $\Phi(T)=(T_n,H_n)$.
While it is not too difficult to show that $T_n$ is a standard Young tableau and that $H_n$ is a hook tableau, it takes considerably more effort to prove that $\Phi$ is a bijection.
For details we refer to~\cite{Kratt1999,NPS1997,Sagan2001}.

\smallskip
Given a tableau $T$ we denote by $n(T)$ the number of exchanges performed during the application of the NPS algorithm, that is, the number of times step \textbf{E5} is visited during the construction of all tableaux $T_1,\dots,T_n$.
The \emph{average case complexity} of the NPS algorithm is now defined as
\begin{align*}
C(\lambda)=\frac{1}{n!}\sum_{T\in\op{T}(\lambda)}n(T).
\end{align*}
The \emph{worst case complexity} is defined as
\begin{align*}
W(\lambda)=\max_{T\in\op{T}(\lambda)} n(T).
\end{align*}

\smallskip
We remark that due to the chosen order of the cells of $\lambda$ the defined algorithm is also called the column-wise NPS algorithm.
In principle, other orders can be chosen and some will result in bijections.
Another result anticipated by Krattenthaler and M{\"u}ller is that the average case complexities of the column-wise and the row-wise NPS algorithms on a fixed shape $\lambda$ are the same.
In other words $C(\lambda)=C(\lambda')$.
See~\cite{NeuSul2015} for a proof of this result and more details in this direction.

\section{Convergence and hook coordinates}\label{Section:hookcoord}

\begin{figure}[t]
\begin{center}
\begin{tikzpicture}[scale=.4]
\draw[gray]
(-1,1)--(3,5)
(-2,2)--(0,4)
(-3,3)--(-2,4)
(-4,4)--(-3,5)
(-5,5)--(-4,6)
(1,1)--(-2,4)
(2,2)--(0,4)
(3,3)--(1,5);
\draw[thick] (0,0)--(4,4)--(2,6)--(0,4)--(-1,5)--(-2,4)--(-5,7)--(-6,6)--cycle;
\draw[thick,->](0,0)--(0,8)node[anchor=west]{$y$};
\draw[thick,->](0,0)--(-8,8)node[anchor=east]{$v$};
\draw[thick,->](0,0)--(8,8)node[anchor=west]{$u$};
\draw[thick,->](-10,0)--(10,0)node[anchor=west]{$x$};
\draw[thick,<->,blue](-8,8)--(-6,6)--(-5,7)--(-2,4)--(-1,5)--(0,4)--(2,6)--(4,4)--(8,8);
\end{tikzpicture}
\caption{The partition $\lambda=(4,4,2,1,1,1)$ depicted in Russian convention.}
\end{center}
\end{figure}

We shall work with the following coordinate system corresponding to a rotation by $\pi/2$ and a rescaling by $\sqrt{n}$
\begin{align}\label{eq:uv}
u = u(x,y) = \frac{\sqrt{n}}{\sqrt{2}}(x+y), \qquad
v = v(x,y) = \frac{\sqrt{n}}{\sqrt{2}}(y-x), \qquad
\left|\frac{\partial(u,v)}{\partial(x,y)}\right| = n.
\end{align}
Given a partition $\lambda$ of $n$ define
\begin{align}\label{eq:Dlambda}
D_{\lambda}
&=\big\{(x,y)\in\R^2:0<v,0<u\leq\lambda_{\ceil{v}}\text{ or }0\leq u,0\leq v,uv=0\big\}.
\end{align}
We define the \emph{boundary function} $\gamma:\R\to\R$ of $\lambda$ via
\begin{align}\label{eq:gamma}
\gamma(x)&=\sup\{y\in\R:(x,y)\in D_{\lambda}\}.
\end{align}
The set $D_{\lambda}$ and function $\gamma$ describe the (rescaled) partition $\lambda$ in the so-called Russian convention.

\smallskip
Let $(\lambda^{(n)})_{n\in\N}$ be a sequence of partitions, such that $\lambda^{(n)}$ is a partition of $n$ and has the boundary function $\gamma_n$.
Denote by $\Gamma$ the set of all $1$-Lipschitz functions $\gamma:\R\to\R$ such that there exists an interval $(a,b)$ with $\gamma(x)=\abs{x}$ for all $x\notin(a,b)$.
Moreover, denote by $\Gamma_1$ the set of all functions $\gamma\in\Gamma$ such that
\begin{align*}
\int_{\R}\gamma(x)-\abs{x}\d x=1.
\end{align*}
Clearly $\gamma_n\in\Gamma_1$ for all $n\in\N$.
We say the sequence $\lambda^{(n)}$ \emph{converges uniformly} to a limit curve $\gamma\in\Gamma_1$ if
\begin{align*}
\lim_{n\to\infty}\sup_{x\in\R}\abs{\gamma(x)-\gamma_n(x)}=0
\end{align*}
and there exists a uniform interval $(a,b)$ such that $\gamma_n(x)=\abs{x}$ for all $x\notin(a,b)$ and all $n\in\N$.

\smallskip
For any $\gamma\in\Gamma$ set
\begin{align*}
D_{\gamma}
&=\big\{(x,y)\in\R^2:\abs{x}\leq y\leq\gamma(x)\big\},
\end{align*}
and let $(x,y)$ be an interior point of $D_{\gamma}$.
We define three functions via
\begin{align*}
a_{\gamma}(x,y)
&=\sqrt{2}\sup\{t\in\R:(x+t,y+t)\in D_{\gamma}\},\\
\ell_{\gamma}(x,y)
&=\sqrt{2}\sup\{t\in\R:(x-t,y+t)\in D_{\gamma}\},\\
d_{\gamma}(x,y)
&=\sqrt{2}\sup\{s+t:s,t>0,(x+s-t,y+s+t)\in D_{\gamma}\}.
\end{align*}
Geometrically $a_{\gamma}(x,y)$ is the distance from $(x,y)$ to the curve $\gamma$ in $u$-direction, $\ell_{\gamma}(x,y)$ is the distance from $(x,y)$ to the curve $\gamma$ in $v$-direction, and $d_{\gamma}(x,y)$ is half of the maximal perimeter among all rectangles confined in $D_{\gamma}$ with sides parallel to the $u$- and $v$-axes, whose lower corner is $(x,y)$.
Extend $a_{\gamma}, \ell_{\gamma}$ and $d_{\gamma}$ from the interior of $D_{\gamma}$ to $\R^2$ by assigning zero to all other points.

\smallskip
With regard to the subsequent sections we need to give the following question some thought.
Suppose $\gamma,\eta\in\Gamma$ are close with respect to the supremum norm, then what can be said about the relationship between $a_{\gamma}$ and $a_{\eta}$ on their common domain $D_{\gamma} \cap D_{\eta}$?
The following example demonstrates that $||a_{\gamma}-a_{\eta}||_{\infty}$ does not need to be small.
Let
\begin{align*}
\gamma(x)
= \begin{cases} \sqrt{2}n+x &\quad\text{if }-\frac{n}{\sqrt{2}}\leq x \leq-\frac{1}{\sqrt{2}}\left(n-\frac{1}{n}\right), \\
\frac{\sqrt{2}}{n}-x &\quad\text{if } -\frac{1}{\sqrt{2}}\left(n-\frac{1}{n}\right)\leq x\leq\frac{1}{n\sqrt{2}}, \\
|x| &\quad\text{else},
\end{cases}\qquad\text{and}\quad\quad
\eta(x)=\gamma(-x).
\end{align*}
Then $||\gamma-\eta||_{\infty} = \frac{\sqrt{2}}{n}$ but $||a_{\gamma}-a_{\eta}||_{\infty} = ||\ell_{\gamma}-\ell_{\eta}||_{\infty} = n-\frac{1}{n}$.

We show, however, in \refl{arm} that when $\gamma$ and $\eta$ agree outside of a fixed interval $(a,b)$, then the exceptional set of points on which $a_{\gamma}$ and $a_{\eta}$ diverge is small when $||\gamma-\eta||_{\infty}$ is small.
The proof of \refl{arm} is based on a geometric argument.
In \refl{leg} and \refl{dist} we deduce analogous results for $\ell_{\gamma}$ and $d_{\gamma}$, which causes little effort once \refl{arm} is established.

For a (measurable) subset $A$ of $\R^n$, let $|A|$ denote the Lebesgue measure of $A$.

\begin{mylem}{symmdiff} Let $\gamma,\eta\in\Gamma$ such that $\gamma(x) = |x| = \eta(x)$ for all $x\notin(a,b)$. Then the Lebesgue measure of the symmetric difference of the sets $D_{\gamma}$ and $D_{\eta}$ is bounded by
\begin{align*}
\big|(D_{\gamma}-D_{\eta})\cup (D_{\eta}-D_{\gamma})\big|
< (b-a) ||\gamma-\eta||_{\infty}.
\end{align*}
\end{mylem}

\begin{proof} It follows immediately from the definition of $D_{\gamma}$ that
\begin{align*}
\big|(D_{\gamma}-D_{\eta})\cup (D_{\eta}-D_{\gamma})\big|
= \int_{a}^{b} |\gamma(x)-\eta(x)|\d x \leq (b-a) ||\gamma-\eta||_{\infty}.
\end{align*}
The inequality is strict since $\gamma$ and $\eta$ agree at $a$ and $b$ and are continuous.
\end{proof}

\begin{mylem}{arm} Let $\varepsilon>0$ and $(a,b)$ be an interval.
Then there exists a constant $K$ such that for all functions $\gamma,\eta\in\Gamma$ with $\gamma(x)=|x|=\eta(x)$ for all $x\notin(a,b)$ we have
\begin{align*}
\big|\big\{(x,y) \in D_{\gamma}\cap D_{\eta} : |a_{\gamma}(x,y)-a_{\eta}(x,y)| > \varepsilon \big\}\big|
\leq \frac{K}{\varepsilon} ||\gamma-\eta||_{\infty}.
\end{align*}
\end{mylem}

\begin{proof} Let $||\gamma-\eta||_{\infty} = \delta$.

We first prove the claim under the assumption that $\gamma(x) \leq \eta(x)$ for all $x\in(a,b)$.
To this end let $k=\ceil{\frac{b-a}{\delta}}$ and $x_0, x_1, \dots, x_k$ be a subdivision of the interval $[a,b]$ such that $x_0=a$, $x_k=b$ and $x_i-x_{i-1}=\frac{b-a}{k}\leq\delta$.
Furthermore, subdivide the curve of $\gamma$ into segments $\sigma_1,\dots,\sigma_k$ such that
\begin{align*}
\sigma_i = \big\{(x,\gamma(x)) : x\in[x_{i-1},x_{i}]\big\}
\qquad\text{for }i\in\{1,\dots,k\}.
\end{align*}
Define
\begin{align*}
\tau_i = \big\{(x,\gamma(x)+\delta) : x\in[x_{i-1},x_{i}]\big\}
\qquad\text{for }i\in\{1\dots,k\},
\end{align*}
and add $\tau_0= \{(x,a+\delta) : x\in[a-\delta,a]\}$ and $\tau_{k+1} = \{(x,b+\delta) : x\in[b,b+\delta]\}$.

A crucial observation is that $\gamma(x)\leq \eta(x)\leq \tau(x)$ for all $x \in \R$ and consequently $a_{\gamma}(x,y) \leq a_{\eta}(x,y) \leq a_{\tau}(x,y)$ for all $(x,y)\in D_{\gamma}$.
Here $\tau\in\Gamma$ denotes the concatenation of the segments $\tau_i$ extended by $\tau(x)=\abs{x}$ for all $x\notin(a-\delta,b+\delta)$.

Let $P:\R^2\to\R^2$ be the orthogonal projection onto the line $\{(x,-x):x\in\R\}$.
Set $P(\sigma_i) = \{P(x,y):(x,y)\in\sigma_{i}\}$ and define $P(\tau_{i})$ analogously.
If $(x,y),(x',y')$ are points in $D_{\gamma}$ such that $P(x,y)=P(x',y')$ then it follows that $a_{\gamma}(x,y)-a_{\eta}(x,y)=a_{\gamma}(x',y')-a_{\eta}(x',y')$.
Thus we should investigate the set of points $\{(x,|x|):a\leq x\leq 0\}$.
Note that the sets $P(\sigma_i)$ cover the set $\{(x,|x|) : a\leq x\leq 0\}$ in such a way that $P(\sigma_i)\cap P(\sigma_{i+1})$ consists of a single point, and $P(\sigma_i)\cap P(\sigma_j)$ is empty unless all sets $P(\sigma_{i+1}),\dots,P(\sigma_{j-1})$ collapse to a single point.

Suppose that $(x,|x|)\in P(\sigma_i)$ with $|a_{\gamma}(x,|x|)-a_{\eta}(x,|x|)| > \varepsilon$ then $(x,|x|)\in P(\sigma_i)\cap P(\tau_{i+m})$ for some $m$ with
\begin{align}\label{eq:m}
m \sqrt{2}\delta
&\geq\varepsilon.
\end{align}
In this case observe that 
\begin{align}\label{eq:Psigma}
\Bigg|\bigcup_{j=0}^{m} P(\sigma_{i+j}) \Bigg|
\leq\frac{5}{\sqrt{2}}\delta.
\end{align}
Choose sequences $(i_1,i_2,\dots)$ and $(m_1,m_2,\dots)$ as follows.
Let $i_1$ be the minimal $i\in\{1,\dots,k\}$ such that for some $(x,|x|) \in P(\sigma_i)$ we have $|a_{\gamma}(x,|x|)-a_{\eta}(x,|x|)| > \varepsilon$.
Given $i_j$ let $m_j$ be maximal such that $P(\sigma_{i_j})\cap P(\tau_{i_{j}+m_j}) \neq \emptyset$.
Given $i_j$ and $m_j$, if there is an $i\in\{i_{j}+m_{j}+1,\dots, k\}$ such that $|a_{\gamma}(x,|x|)-a_{\eta}(x,|x|)| > \varepsilon$ for some $(x,|x|)\in P(\sigma_i)$, then let $i_{j+1}$ be the minimal $i$ with this property.
Otherwise terminate both sequences.

Clearly the sequence $(i_1,i_2,\dots,i_r)$ is finite.
Because of \refq{m} we have $m_j \geq \frac{\varepsilon}{\delta\sqrt{2}}$ for all $j$, and hence
\begin{align*}
(j-1)\frac{\varepsilon}{\delta\sqrt{2}}
&<i_j
\leq k
<\frac{b-a+1}{\delta}.
\end{align*}
In particular,
\begin{align}\label{eq:r}
r
&<\frac{\sqrt{2}(b-a+1)}{\varepsilon}+1.
\end{align}
By definition of the sequences $(i_j)_j$ and $(m_j)_j$ for every point $(x,y) \in D_{\gamma}$ with $|a_{\gamma}(x,y)-a_{\eta}(x,y)|>\varepsilon$ the projection $P(x,y)$ is contained in $P(\sigma_{i_j})\cup\dots\cup P(\sigma_{i_j+m_j})$ for some $j\in\{1,\dots,r\}$.
But now, using \refq{Psigma} and \refq{r} we obtain
\begin{align}\notag
&\big|\big\{ (x,y) : |a_{\gamma}(x,y)-a_{\eta}(x,y)|>\varepsilon \big\}\big|\\
\notag
&\qquad\leq \big|\big\{ (x,y) : P(x,y)\in P(\sigma_{i_j})\cup\dots\cup P(\sigma_{i_{j}+m_{j}})\text{ for some }1\leq j \leq r \big\}\big|\\
\notag
&\qquad< rb\sqrt{2}\frac{5}{\sqrt{2}}\delta\\
\label{eq:leg}
&\qquad< \frac{1+5b(b-a+1)\sqrt{2}}{\varepsilon}\delta.
\end{align}
Now drop the condition $\gamma\leq \eta$, and assume $||\gamma-\eta||_{\infty} < \frac{\delta}{2}$.
Consider the function $\rho$ defined by $\rho(x)=\max\{\gamma(x)-\frac{\delta}{2},|x|\}$.
Clearly $\rho\leq\gamma$, $||\rho-\gamma||_{\infty}<\delta$ and $\rho\leq\eta$, $||\rho-\eta||_{\infty}<\delta$. Thus, appealing to \refq{leg} twice, the inclusion 
\begin{align*}
&\big\{ (x,y) : |a_{\gamma}(x,y)-a_{\eta}(x,y)|>\varepsilon \big\} \\
&\qquad\subseteq \big\{ (x,y) : |a_{\rho}(x,y)-a_{\gamma}(x,y)|>\varepsilon \big\}
\cup \big\{ (x,y) : |a_{\rho}(x,y)-a_{\eta}(x,y)|>\varepsilon \big\} \\
&\qquad\quad\cup (D_{\rho}-D_{\gamma})\cup (D_{\gamma}-D_{\rho})
\cup (D_{\rho}-D_{\eta})\cup (D_{\eta}-D_{\rho})
\end{align*}
holds and \refl{symmdiff} completes the proof.
\end{proof}

The analogous result on $\ell_{\gamma}$ follows easily.
\begin{mylem}{leg} Let $\varepsilon > 0$ and $(a,b)$ be an interval.
Then there exists a constant $K$ such that for all $\gamma,\eta\in\Gamma$ with $\gamma(x)=|x|=\eta(x)$ for all $x\notin(a,b)$ we have
\begin{align*}
\big|\big\{ (x,y)\in D_{\gamma}\cap D_{\eta} : |\ell_{\gamma}(x,y)-\ell_{\eta}(x,y)| > \varepsilon \big\}\big|
< \frac{K}{\varepsilon} ||\gamma-\eta||_{\infty}.
\end{align*}
\end{mylem}

\begin{proof} The claim follows directly from \refl{arm} and the symmetry
\begin{align*}
a_{\gamma}(x,y) = \ell_{\rho}(-x,y),
\end{align*}
where $\rho\in\Gamma$ is the function defined by $\rho(x) = \gamma(-x)$.
\end{proof}

Finally there is a similar result for the function $d_{\gamma}$.
\begin{mylem}{dist} Let $\varepsilon>0$ and $(a,b)$ be an interval.
Then there exists a constant $K$ such that for all $\gamma,\eta\in\Gamma$ with $\gamma(x)=|x|=\eta(x)$ for all $x\notin(a,b)$ we have
\begin{align*}
\big|\big\{ (x,y) \in D_{\gamma}\cap D_{\eta} : |d_{\gamma}(x,y)-d_{\eta}(x,y)| > \varepsilon \big\}\big|
< \frac{K}{\varepsilon} ||\gamma-\eta||_{\infty}.
\end{align*}
\end{mylem}

\begin{proof} The claim follows from \refl{arm} and \refl{leg} and the estimation
\begin{align*}
|d_{\gamma}(x,y)-d_{\eta}(x,y)|
&\leq |a_{\gamma}(x,y)-a_{\eta}(x,y)| + |\ell_{\gamma}(x,y)-\ell_{\eta}(x,y)| + \sqrt{2} ||\gamma-\eta||_{\infty}.
\end{align*}
\end{proof}

\begin{figure}[ht]
\begin{center}
\begin{tikzpicture}[scale=.4]
\draw[thick,dashed]
(1,2)--(-1.5,4.5)--(-1.5,0)node[anchor=north]{$s$}
(1,2)--(3.5,4.5)--(3.5,0)node[anchor=north]{$t$};
\draw[thick,->](0,0)--(0,8);
\draw[thick,->](0,0)--(-8,8);
\draw[thick,->](0,0)--(8,8);
\draw[thick,->](-10,0)--(10,0);
\draw[thick,<->,blue](-8,8)--(-6,6)--(-5,7)--(-2,4)--(-1,5)--(0,4)--(2,6)--(4,4)--(8,8);
\draw (1,2)node{$\bullet$};
\end{tikzpicture}
\caption{The hook coordinates $(s,t)$ of an interior point $(x,y)\in D_{\lambda}$.}
\end{center}
\end{figure}
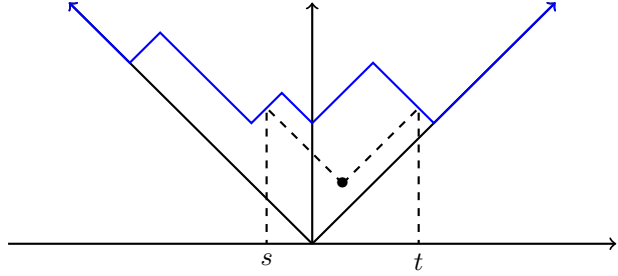

\smallskip
We conclude this section by introducing the so-called hook coordinates, which were named (to the best of the authors' knowledge) by Dan Romik, and appear naturally in the study of limit shapes of partitions, see for example \cite{LogShe1977}. Namely, we set
\begin{align}\label{eq:hookcoord}
s=x-\frac{\ell_{\gamma}(x,y)}{\sqrt{2}},\qquad
t=x+\frac{a_{\gamma}(x,y)}{\sqrt{2}}.
\end{align}
Note that
\begin{align*}
s+\gamma(s)=x+y,\qquad
t-\gamma(t)=x-y,
\end{align*}
which yields
\begin{align*}
x=\frac{1}{2}(s+t+\gamma(s)-\gamma(t)),\qquad
y=\frac{1}{2}(s-t+\gamma(s)+\gamma(t)),
\end{align*}
and
\begin{align*}
\left|\frac{\partial(x,y)}{\partial(s,t)}\right|
&=\frac{1}{2}(1+\gamma'(s))(1-\gamma'(t)),
\end{align*}
where the derivative $\gamma'$ is defined almost everywhere since $\gamma$ is 1-Lipschitz.

\section{Worst case complexity}\label{Section:W}

In this section we analyze the asymptotic behavior of the worst case complexity of the NPS algorithm.
We first demonstrate in \refp{W} that a trivial combinatorial upper bound for the worst case complexity of the NPS algorithm on a fixed shape $W(\lambda)$ is in fact tight.
Furthermore, \reft{Wn} provides the first order asymptotics of $W(\lambda^{(n)})$, where $(\lambda^{(n)})_{n\in\N}$ converges uniformly, in terms of the limit curve $\gamma$.

\smallskip
For a cell $(i,j)\in\lambda$ denote by
\begin{align}\label{eq:W}
w(i,j)
=\max\big\{i'-i+j'-j:(i',j')\in\lambda,i\leq i'\text{ and }j\leq j'\big\}
\end{align}
the maximal distance of the cell $(i,j)$ to a cell $(i',j')\in\lambda$ South-East of $(i,j)$.
Let $W(\lambda)$ denote the worst case complexity of the NPS algorithm on $\lambda$.
Clearly
\begin{align*}
W(\lambda)
\leq\sum_{(i,j)\in\lambda}w(i,j).
\end{align*}
We first show that this upper bound is tight.

\begin{figure}[t]
\begin{center}
\begin{tikzpicture}[scale=.5]
\begin{scope}
\draw[gray]
(0,1)--(1,1)
(1,2)--(1,4)
(2,2)--(2,4)
(3,2)--(3,4)
(4,3)--(4,4);
\draw[very thick]
(0,0)rectangle(1,2)
(0,2)rectangle(4,3)
(0,3)rectangle(5,4);
\draw[xshift=5mm,yshift=5mm]
(4,3)node{{$1$}}
(3,3)node{{$2$}}
(2,3)node{{$3$}}
(1,3)node{{$4$}}
(0,3)node{{$5$}}
(3,2)node{{$6$}}
(2,2)node{{$7$}}
(1,2)node{{$8$}}
(0,2)node{{$9$}}
(0,0)node{{$10$}}
(0,1)node{{$11$}};
\end{scope}
\begin{scope}[xshift=8cm]
\draw[gray]
(0,1)--(2,1)
(0,2)--(2,2)
(1,0)--(1,4)
(2,2)--(2,4)
(3,2)--(3,4)
(4,3)--(4,4);
\draw[very thick]
(0,0)rectangle(2,3)
(2,2)rectangle(4,3)
(0,3)rectangle(5,4);
\draw[xshift=5mm,yshift=5mm]
(4,3)node{{$1$}}
(3,3)node{{$2$}}
(2,3)node{{$3$}}
(1,3)node{{$4$}}
(0,3)node{{$5$}}
(1,0)node{{$6$}}
(1,1)node{{$7$}}
(1,2)node{{$8$}}
(0,0)node{{$9$}}
(0,1)node{{$10$}}
(0,2)node{{$11$}}
(3,2)node{{$12$}}
(2,2)node{{$13$}}
;
\end{scope}
\begin{scope}[xshift=16cm,yshift=-1cm]
\draw[gray]
(0,2)--(4,2)
(0,3)--(4,3)
(0,4)--(4,4)
(1,1)--(1,5)
(2,1)--(2,5)
(3,1)--(3,5);
\draw[very thick]
(0,0)rectangle(1,1)
(0,1)rectangle(4,5)
(4,4)rectangle(5,5);
\draw[xshift=5mm,yshift=5mm]
(3,1)node{{$1$}}
(3,2)node{{$2$}}
(3,3)node{{$3$}}
(3,4)node{{$4$}}
(2,1)node{{$5$}}
(2,2)node{{$6$}}
(2,3)node{{$7$}}
(2,4)node{{$8$}}
(1,1)node{{$9$}}
(1,2)node{{$10$}}
(1,3)node{{$11$}}
(1,4)node{{$12$}}
(0,1)node{{$13$}}
(0,2)node{{$14$}}
(0,3)node{{$15$}}
(0,4)node{{$16$}}
(4,4)node{{$17$}}
(0,0)node{{$18$}};
\end{scope}
\end{tikzpicture}
\caption{Three tableaux that exhibit the worst case complexity of the NPS algorithm on their respective shapes.}
\label{Figure:worst}
\end{center}
\end{figure}

\begin{myprop}{W}
Let $\lambda$ be a partition. Then
\begin{align}\label{eq:Wlambda}
W(\lambda)=\sum_{(i,j)\in\lambda}w(i,j)
\end{align}
\end{myprop}

\begin{proof}
We construct an explicit tableau $T\in\op{T}(\lambda)$ such that the number of exchanges $n(T)$ equals the right hand side of \refq{Wlambda}.

If $\lambda$ is a rectangle, that is, $\lambda$ has only one corner, then let $(i_1,j_1)\succ\dots\succ(i_n,j_n)$ be the cells of $\lambda$.
Define the tableau $T$ by setting $T(i_r,j_r)=r$.
During the application of the NPS algorithm to $T$ every entry is moved to the corner of $\lambda$ and \refq{Wlambda} holds.

\smallskip
If $\lambda$ is of general form we construct $T$ and a sequence of cells $(i_1,j_1),\dots,(i_l,j_l)$ as follows.
Set $(i_1,j_1)=(1,1)$.
Given $(i_r,j_r)$ let
\begin{align*}
X_r=\big\{(i,j)\in\lambda:i_r\leq i,j_r\leq j,w(i_r,j_r)=i-i_r+j-j_r\big\}
\end{align*}
denote the set of corners South-East of $(i_r,j_r)$ with maximal distance to $(i_r,j_r)$.
Fix any corner $(i_r',j_r')\in X_r$ and let
\begin{align*}
R_r=\big\{(i,j)\in\lambda:i_r\leq i\leq i_r',j_r\leq j\leq j_r'\big\}
\end{align*}
denote the rectangle inside $\lambda$ that contains $(i_r,j_r)$ and $(i_r',j_r')$.
Now define $T$ on $R_r$ by assigning the numbers 
\begin{align*}
\Big\{\Big(\sum_{k=1}^{r-1}(i_k'-i_k+1)(j_k'-j_k+1)\Big)+1,\dots,\sum_{k=1}^{r}(i_k'-i_k+1)(j_k'-j_k+1)\Big\}
\end{align*}
to the cells in the rectangle $R_r$ using the reverse lexicographic order as above.
If $T$ is not defined on all cells of $\lambda$ then let $(i_{r+1},j_{r+1})$ be a cell of $\lambda$ with maximal hook-length among all cells for which $T$ is not yet defined.
Compare to \reff{worst}.

By construction each entry $T(i,j)$ of a cell $(i,j)\in R_r$ drops to the corner $(i_r',j_r')$ of the rectangle $R_r$ during the application of the NPS algorithm.
Thus $n(T)$ is given by the right hand side of \refq{Wlambda}.
\end{proof}

Approximating the right hand side of~\refq{W} by an integral and making use of the preparatory results in \refs{hookcoord} we are able to draw conclusions on the asymptotics of the worst case complexity.

\begin{mythm}{Wn} Let $(\lambda^{(n)})_{n\in\N}$ be a sequence of partitions converging uniformly to the limit shape $\gamma\in\Gamma_1$. Then
\begin{align}\label{eq:Wn}
W(\lambda^{(n)})
&=n^{3/2}\iint_{D_{\gamma}}d_{\gamma}(x,y)\d x\d y+o(n^{3/2})&&\text{as }n\to\infty.
\end{align}
\end{mythm}

Before we turn to the proof let us state some remarks.

First, let us argue the existence of the integral in \refq{Wn}.
Since $d_{\gamma}(x,y)$ is bounded and $D_{\gamma}$ is the union of a compact set and a null set, the integral is proper.
Furthermore, the function $d_{\gamma}(x,y)$ is decreasing in $y$ and hence integrable.
The function $\int_{\abs{x}}^{\gamma(x)}d_{\gamma}(x,y)\d y$ is even continuous in $x$.

Secondly, since the right hand side of \refq{Wn} is \emph{a priori} not straight forward to compute, we offer the estimation
\begin{align*}
\iint_{D_{\gamma}}d_{\gamma}(x,y)\d x\d y
&\leq \frac{\sqrt{2}}{2}\int_{-\infty}^{\infty}\int_{s}^{\infty} (t-s)\big(1+\gamma'(s)\big)\big(1-\gamma'(t)\big)\d t\d s,
\end{align*}
which is obtained from $d_{\gamma}(x,y)\leq a_{\gamma}(x,y)+\ell_{\gamma}(x,y)$ by a substitution of hook coordinates.

\begin{proof}[Proof of \reft{Wn}] First rewrite $W(\lambda^{(n)})$ as an integral as follows:
A cell $(i,j)\in\lambda^{(n)}$ corresponds to the square
\begin{align*}
Z(i,j)
=\big\{(x,y):i-1\leq v\leq i,j-1\leq u\leq j\big\}
\subseteq D_n.
\end{align*}
Let $(x,y)$ be an interior point of $Z(i,j)$.
Then
\begin{align*}
\sqrt{n}d_{\gamma_n}(x,y)
&=w(i,j)+i+j-u-v.
\end{align*}
Hence
\begin{align*}
n^{3/2}\iint_{Z(i,j)}d_{\gamma_n}(x,y)\d x\d y
&=n\iint_{Z(i,j)}w(i,j)+i+j-u-v\d x\d y\\
&=\int_{i-1}^i\int_{j-1}^{j}w(i,j)+i+j-u-v\d u\d v\\
&=w(i,j)+1
\end{align*}
and
\begin{align*}
\sum_{(i,j)\in\lambda^{(n)}}w(i,j)
&=-n+n^{3/2}\iint_{D_n}d_{\gamma_n}(x,y)\d x\d y.
\end{align*}
Now fix $\varepsilon>0$.
It suffices to show that
\begin{align*}
\abs{\iint_{D_{\gamma}}d_{\gamma}(x,y)\d x\d y
-\iint_{D_n}d_{\gamma_n}\d x\d y}
&<\varepsilon
\end{align*}
for all sufficiently large $n$.
In order to do so choose an interval $(a,b)$ such that $\gamma_n(x)=\abs{x}$ for all $x\notin(a,b)$ and all $n\in\N$.
It follows that also $\gamma(x)=|x|$ outside of $(a,b)$.
By \refl{symmdiff} the Lebesgue measure of the symmetric difference of the sets $D_{\gamma}$ and $D_n$ tends to zero as $n$ tends to infinity.
Since both $d_{\gamma}$ and $d_{\gamma_n}$ are bounded by the constant $(b-a)\sqrt{2}$, 
\begin{align*}
\iint_{D_n-D_{\gamma}}d_{\gamma_n}(x,y)\d x\d y
+\iint_{D_{\gamma}-D_n}d_{\gamma}(x,y)\d x\d y
&<\varepsilon
\end{align*}
when $||\gamma-\gamma_n||_{\infty}<\varepsilon/((b-a)^2\sqrt{2})$.
On the other hand by \refl{dist} there exist sets $A$ and $B$ and a constant $K$ such that $D_{\gamma}\cap D_{n}=A\cup B$, $|d_{\gamma}(x,y)-d_{\gamma_n}(x,y)|<\varepsilon/2$ for all $(x,y)\in A$ and $|B|<K||\gamma-\gamma_n||_{\infty}/\varepsilon$.
Hence,
\begin{align*}
\iint_{D_{\gamma}\cap D_n}\abs{d_{\gamma}(x,y)-d_{\gamma_n}(x,y)}\d y\d x
&<\frac{\varepsilon}{2}+(b-a)\sqrt{2}\frac{K||\gamma-\gamma_n||_{\infty}}{\varepsilon}
<\varepsilon
\end{align*}
if $||\gamma-\gamma_n||_{\infty}$ is sufficiently small. 
\end{proof}

\section{Average case complexity}\label{Section:C}

The main result of this section is an asymptotic lower bound for the average case complexity of the NPS algorithm, which we obtain in three steps.
\refp{EH} obtains a combinatorial bound for the average case complexity of the NPS algorithm on a fixed shape $C(\lambda)$.
\refp{C} approximates this combinatorial bound by an integral.
Finally, in \reft{Cn} we derive an asymptotic bound for $C(\lambda^{(n)})$, where $(\lambda^{(n)})_{n\in\N}$ converges uniformly, in terms of the limit curve $\gamma$.

\smallskip
Given a hook tableau $H$ of shape $\lambda$ denote
\begin{align*}
|H|
&=\sum_{(i,j)\in\lambda}\abs{H(i,j)}.
\end{align*}

\begin{myprop}{EH} Let $\lambda$ be a partition of $n$, and $H$ be a hook tableau of shape $\lambda$ chosen uniformly at random. Then
\begin{equation*}
C(\lambda)>\E(\abs{H}).
\end{equation*}
\end{myprop}

\begin{proof} During the application of the NPS algorithm to a tableau $T$ of shape $\lambda$ a hook tableau of the same shape is built from the zero tableau, that is, the hook tableau with $\abs{H}=0$.
This is done by applying the following transformations.
Suppose the entry of the cell $(i,j)$ drops to the cell $(i',j')$.
Then $H(s,j)$ is set to $H(s+1,j)-1$ for $s=i,\dots,i'-1$, and $H(i',j)$ is set to $j'-j$.
Thereby $|H|$ is increased by no more than $i'-i+j'-j$ which is exactly the number of performed exchanges.
Since the NPS algorithm produces each hook tableau equally often as $T$ ranges over all possible tableaux of shape $\lambda$, we conclude the following lower bound on the average number of exchanges
\begin{align}\label{eq:E}
C(\lambda)>\frac{f^{\lambda}}{n!}\sum_{H}|H|,
\end{align}
where the sum is taken over all hook tableaux of shape $\lambda$.
The number of hook tableaux is given by the hook product
\begin{align*}
\prod_{(i,j)\in\lambda}h_{\lambda}(i,j).
\end{align*}
By use of the hook-length formula, the right hand side of \refq{E} is just the expected value of the random variable $|H|$.
\end{proof}

In a next step we replace the combinatorial lower bound $C(\lambda)>\E(\abs{H})$ by an integral.

\begin{myprop}{C} Let $\lambda$ be a partition of $n$ with boundary $\gamma\in\Gamma_1$.
Then
\begin{align*}
C(\lambda)
&>\frac{n^{3/2}}{2}\iint_{D_{\gamma}}\frac{a_{\gamma}(x,y)^2+\ell_{\gamma}(x,y)^2}{a_{\gamma}(x,y)+\ell_{\gamma}(x,y)+\frac{1}{\sqrt{n}}}\d x\d y+\frac{n}{2}.
\end{align*}
\end{myprop}

\begin{proof} Recall that $C(\lambda)>\E(\abs{H})$, where $H$ ranges over the hook tableaux of shape $\lambda$ by \refp{EH}.
By linearity
\begin{align*}
\E(|H|)
&=\sum_{(i,j)\in\lambda}\E(\abs{H(i,j)})\\
%
%
&=\sum_{(i,j)\in\lambda}\frac{\arm(i,j)^2+\arm(i,j)+\leg(i,j)^2+\leg(i,j)}{2h(i,j)}\\
&=\bigg(\sum_{(i,j)\in\lambda}\frac{\arm(i,j)^2+\leg(i,j)^2}{2h(i,j)}\bigg)+\frac{n}{2}-\sum_{(i,j)\in\lambda}\frac{1}{2h(i,j)}.
\end{align*}
For any cell $(i,j)\in\lambda$ let $Z(i,j)=\{(x,y):i-1\leq v\leq i,j-1\leq u\leq j\}$ as in the proof of \reft{Wn}.
If $(x,y)$ is an interior point of $Z(i,j)$ then
\begin{align*}
a_{\gamma}(x,y)\sqrt{n}
&=\arm(i,j)+j-u
&&\text{and}&&
\ell_{\gamma}(x,y)\sqrt{n}=
\leg(i,j)+i-v.
\end{align*}
Throughout the remainder of this proof denote $a_{ij}=\arm(i,j)$, $l_{ij}=\leg(i,j)$, $h_{ij}=h(i,j)$, $w=i-v$ and $z=j-u$.
A straightforward computation yields
\begin{align*}
&\frac{n}{\sqrt{n}}\frac{a_{\gamma}(x,y)^2}{a_{\gamma}(x,y)+\ell_{\gamma}(x,y)+\frac{1}{\sqrt{n}}} - \frac{a_{ij}^2}{h_{ij}}
= \frac{(a_{ij}+z)^2}{h_{ij}+w+z} - \frac{a_{ij}^2}{h_{ij}} \\
&\qquad=\frac{(a_{ij}+l_{ij}+1)(a_{ij}^2+2a_{ij}z+z^2)-(a_{ij}+l_{ij}+w+z+1)a_{ij}^2}{h_{ij}(h_{ij}+w+z)} \\
&\qquad=\frac{a_{ij}^2z+a_{ij}z^2+2a_{ij}l_{ij}z+l_{ij}z^2+2a_{ij}z+z^2-a_{ij}^2w}{h_{ij}(h_{ij}+w+z)} \\
&\qquad=\frac{a_{ij}^2}{h_{ij}} \frac{z-w}{h_{ij}+w+z} + \frac{a_{ij}l_{ij}+a_{ij}}{h_{ij}}\frac{2z}{h_{ij}+w+z} + \frac{a_{ij}+l_{ij}+1}{h_{ij}}\frac{z^2}{h_{ij}+w+z}.
\end{align*}
Thus
\begin{align}
\frac{a_{ij}^2}{h_{ij}}
&=n\iint_{Z(i,j)}\frac{a_{ij}^2}{h_{ij}}\d x\d y\notag\\
&=n^{3/2}\iint_{Z(i,j)}\frac{a_{\gamma}(x,y)^2}{a_{\gamma}(x,y)+\ell_{\gamma}(x,y) +\frac{1}{\sqrt{n}}}\d x\d y\notag\\
\label{eq:err1}
&\qquad-n\frac{a_{ij}^2}{h_{ij}}\iint_{Z(i,j)}\frac{w-z}{h_{ij}+w+z}\d x\d y\\
\label{eq:err2}
&\qquad-n\frac{a_{ij}l_{ij}+a_{ij}}{h_{ij}^2}\iint_{Z(i,j)}\frac{2h_{ij}z}{h_{ij}+w+z}\d x\d y\\
\label{eq:err3}
&\qquad-n\iint_{Z(i,j)}\frac{z^2}{h_{ij}+w+z}\d x\d y.
\end{align}
The error term \refq{err1} vanishes by symmetry in $z$ and $w$. We have
\begin{align*}
n\iint_{Z(i,j)}\frac{w-z}{h_{ij}+w+z}\d y\d x
&=\int_{0}^{1}\int_{0}^{1}\frac{w-z}{h_{ij}+w+z}\d z\d w
=0.
\end{align*}
The other two error terms, \refq{err2} and \refq{err3}, could be neglected as the integrands are non-negative.
However, there is no harm in showing that they are also small. 
This can be seen easily since the integrands converge uniformly to $2z$ and $0$ respectively.
Thus
\begin{align*}
\lim_{m\to\infty} \int_{0}^{1}\int_{0}^{1} \frac{2mz}{m+w+z} \d z\d w = 1,
\qquad
\lim_{m\to\infty} \int_{0}^{1}\int_{0}^{1} \frac{z^2}{m+w+z}\d z\d w = 0.
\end{align*}
In particular, the integrals in \refq{err2} and \refq{err3} are uniformly bounded no matter how large $h_{ij}$ is.

Approximating $\ell_{ij}^2/h_{ij}$ by an analogous integral and summing over all cells of $\lambda$ yields the claim.
\end{proof}

Using the preparatory results of \refs{hookcoord}, we obtain the main theorem of this section.

\begin{mythm}{Cn} Let $(\lambda^{(n)})_{n\in\N}$ be a sequence of partitions converging uniformly to the limit curve $\gamma\in\Gamma_1$.
Then
\begin{align}\label{eq:Cn}
C(\lambda^{(n)})
&>\frac{\sqrt{2}n^{3/2}}{8}
\int_{-\infty}^{\infty}\int_{s}^{\infty}\left((t-s)+ \frac{(\gamma(t)-\gamma(s))^2}{t-s}\right) \big(1+\gamma'(s)\big)\big(1-\gamma'(t)\big)\d t \d s+o(n^{3/2})
\end{align}
as $n\to\infty$.
\end{mythm}

Before we give a proof let us argue that the right hand side of \refq{Cn} is well-defined.
First note that the integral is really taken over a compact set.
Suppose that $\gamma(x)=\abs{x}$ for all $x\notin(a,b)$ then $\gamma'(s)=-1$ for all $s<a$ and $\gamma'(t)=1$ for all $t>b$, and the integrand vanishes.
Because $\gamma$ is Lipschitz continuous, its derivative exists almost everywhere, is Lebesgue integrable and fulfils $\int_{a}^{b}\gamma'(s)\d s=\gamma(b)-\gamma(a)$.
Thus the limit of the quotient $(\gamma(t)-\gamma(s))/(t-s)$ as $t$ tends to $s$ exists almost everywhere.
The integrand is therefore essentially bounded and integrable.

\begin{proof}[Proof of \reft{Cn}] Let $\gamma_n\in\Gamma_1$ be the boundary function of $\lambda^{(n)}$ and choose an interval $(a,b)$ such that $\gamma_n(x)=|x|$ for all $x\notin(a,b)$.
We begin by noting that
\begin{align}\label{eq:claimv}
\abs{\iint_{D_{\gamma_n}}\frac{a_{\gamma_n}(x,y)^2+\ell_{\gamma_n}(x,y)^2}{a_{\gamma_n}(x,y)+\ell_{\gamma_n}(x,y)+\frac{1}{\sqrt{n}}}\d x\d y
-\iint_{D_{\gamma}}\frac{a_{\gamma}(x,y)^2+\ell_{\gamma}(x,y)^2}{a_{\gamma}(x,y)+\ell_{\gamma}(x,y)}\d x\d y}
\to0
\end{align}
as $n$ tends to infinity.
To see this, first note that the functions $a_{\gamma_n}$, $a_{\gamma}$, $\ell_{\gamma_n}$ and $\ell_{\gamma}$ are all non-negative and bounded by $(b-a)\sqrt{2}$.
Thus also the integrands in \refq{claimv} are non-negative and bounded.
By \refl{symmdiff} it suffices to consider the common domain $D_{\gamma_n}\cap D_{\gamma}$.
It then follows from the Lemmata~\ref{Lemma:arm} and~\ref{Lemma:leg} that
\begin{align}
C(\lambda^{(n)})
&>\frac{n^{3/2}}{2}\iint_{D_{\gamma}}\frac{a_{\gamma}(x,y)^2+\ell_{\gamma}(x,y)^2}{a_{\gamma}(x,y)+\ell_{\gamma}(x,y)}\d x\d y+o(n^{3/2})
&&\text{as }n\to\infty.
\end{align}
To finish the proof we use our hook coordinates.
Substitution gives
\begin{align*}
&\frac{1}{2} \iint_{D_{\gamma}} \frac{a_{\gamma}^2(x,y)+\ell_{\gamma}(x,y)^2}{a_{\gamma}(x,y)+\ell_{\gamma}(x,y)}\d x\d y\\
&\quad= \frac{1}{2} \iint_{D_{\gamma}} \frac{2(t-x)^2+2(x-s)^2}{\sqrt{2}(t-s)}\d x\d y\\
&\quad= \frac{\sqrt{2}}{2} \iint_{D_{\gamma}} \frac{\left(t-\frac{s+t+\gamma(s)-\gamma(t)}{2}\right)^2 + \left(\frac{s+t+\gamma(s)-\gamma(t)}{2}-s\right)^2}{t-s}\d x\d y\\
&\quad= \frac{\sqrt{2}}{8} \iint_{D_{\gamma}} \frac{\big((t-s)+(\gamma(t)-\gamma(s))\big)^2 + \big((t-s)-(\gamma(t)-\gamma(s))\big)^2}{t-s}\d x\d y\\
&\quad= \frac{\sqrt{2}}{8} \iint_{D_{\gamma}} \frac{2(t-s)^2 + 2(\gamma(t)-\gamma(s))^2}{t-s}\d x\d y\\
&\quad= \frac{\sqrt{2}}{8} \int_{-\infty}^{\infty}\int_{s}^{\infty} \left((t-s) + \frac{(\gamma(t)-\gamma(s))^2}{t-s}\right) (1+\gamma'(s))(1-\gamma'(t))\d t \d s.
\end{align*}
Perhaps the only step that needs comment is the last one.
Each pair $(s,t)$ with $s\leq t$ gives rise to a unique point $(x,y)\in D_{\gamma}$ unless $\gamma'(s)=-1$ or $\gamma'(t)=1$.
This follows from the Lipschitz property of $\gamma$.
However, the integrand vanishes when either of the two cases $\gamma'(s)=-1$ or $\gamma'(t)=1$ occurs.
Hence, we can relax the limits of the integral without altering its evaluation.
\end{proof}

\section{Imbalanced scaling}\label{Section:pq}

Some types of partitions, for example partitions with a fixed number of parts, do not converge to a limit curve $\gamma\in\Gamma_1$ in the sense of \refs{hookcoord}.
However, they might converge if an alternative scaling, that is, not by a factor of $\sqrt{n}$ in both $u$- and $v$-direction, is chosen.
In this section we study the asymptotic behaviour of the average case and worst case complexity of the NPS algorithm when the partitions under consideration converge after an imbalanced scaling.

\smallskip
For $p,q\in\Q\cup\{\infty\}$ such that $p,q>0$ and
\begin{align*}
\frac{1}{p}+\frac{1}{q}&=1
\end{align*}
consider the coordinates given by
\begin{align}\label{eq:uvpq}
u=u(x,y)=\frac{n^{1/p}}{\sqrt{2}}(x+y),
&& v=v(x,y)=\frac{n^{1/q}}{\sqrt{2}}(y-x),
&&\abs{\frac{\partial(u,v)}{\partial(x,y)}}=n.
\end{align}
To each partition $\lambda$ of $n$ we associate a set $D_{\lambda}$ and a \emph{$p,q$-boundary function} $\gamma$, defined exactly as in \refq{Dlambda} and \refq{gamma} but with $u$ and $v$ now given by \refq{uvpq}.

Let $(\lambda^{(n)})_{n\in\N}$ be a sequence of partitions such that $\lambda^{(n)}$ is a partition of $n$ and has $p,q$-boundary $\gamma_n$.
We say $(\lambda^{(n)})_{n\in\N}$ converges \emph{$p,q$-uniformly} to the limit curve $\gamma\in\Gamma_1$ if
\begin{align*}
\lim_{n\to\infty}\sup_{x\in\R}\abs{\gamma(x)-\gamma_n(x)}&=0
\end{align*}
and there exists an interval $(a,b)$ such that $\gamma_n(x)=\abs{x}$ for all $x\notin(a,b)$ and all $n\in\N$.

\smallskip
The main result of this section provides the first order asymptotics of $C(\lambda^{(n)})$ and $W(\lambda^{(n)})$, where $(\lambda^{(n)})_{n\in\N}$ converges $p,q$-uniformly, in terms of the limit curve $\gamma$.
They turn out to be of order $n^{1+\max\{1/p,1/q\}}$ as $n$ tends to infinity.

\begin{mythm}{CWpq} Let $(\lambda^{(n)})_{n\in\N}$ be a sequence of partitions converging $p,q$-uniformly to the limit curve $\gamma\in\Gamma_1$.
If $p<q$ then
\begin{align}\label{eq:Cp}
C(\lambda^{(n)})
&=n^{(p+1)/p}\frac{I_1}{2}+o(n^{(p+1)/p})
&&\text{as }n\to\infty,
\end{align}
and
\begin{align}\label{eq:Wp}
W(\lambda^{(n)})
&=n^{(p+1)/p}I_1+o(n^{(p+1)/p})
&&\text{as }n\to\infty,
\end{align}
where
\begin{align*}
I_1
&=\iint_{D_{\gamma}}a_{\gamma}(x,y)\d x\d y
=\frac{\sqrt{2}}{4}\int_{-\infty}^{\infty}\int_{s}^{\infty}(t-s+\gamma(t)-\gamma(s))(1+\gamma'(s))(1-\gamma'(t))\d t\d s.
\end{align*}
On the other hand, if $p>q$ then
\begin{align}\label{eq:Cq}
C(\lambda^{(n)})
&=n^{(p+1)/p}\frac{I_2}{2}+o(n^{(p+1)/p})
&&\text{as }n\to\infty,
\end{align}
and
\begin{align}\label{eq:Wq}
W(\lambda^{(n)})
&=n^{(p+1)/p}I_2+o(n^{(p+1)/p})
&&\text{as }n\to\infty,
\end{align}
where
\begin{align*}
I_2
&=\iint_{D_{\gamma}}\ell_{\gamma}(x,y)\d x\d y
=\frac{\sqrt{2}}{4}\int_{-\infty}^{\infty}\int_{s}^{\infty}(t-s-\gamma(t)+\gamma(s))(1+\gamma'(s))(1-\gamma'(t))\d t\d s.
\end{align*}
\end{mythm}

\begin{proof} Suppose that $p<q$.
Recall from the proof of \refp{C} that
\begin{align}\label{eq:Cpq}
\frac{1}{2}\sum_{(i,j)\in\lambda^{(n)}} \frac{\arm(i,j)^2+\leg(i,j)^2}{\arm(i,j)+\leg(i,j)+1}
+\frac{n}{2}
&<\E(\abs{H})
<C(\lambda^{(n)}).
\end{align}
For any cell $(i,j)\in\lambda^{(n)}$ let $Z(i,j)=\{(x,y):i-1\leq v\leq i,j-1\leq u\leq j\}$.
If $(x,y)$ is an interior point of $Z(i,j)$ then
\begin{align*}
n^{1/p}a_{\gamma_n}(x,y)=\arm(i,j)+j-u
&&\text{and}
&&n^{1/q}\ell_{\gamma_n}(x,y)=\leg(i,j)+i-v,
\end{align*}
where $\gamma_n$ denotes the $p,q$-boundary of $\lambda^{(n)}$.
The asymptotically dominant term in the left hand side of \refq{Cpq} can be approximated by an integral
\begin{align*}
\sum_{(i,j)\in\lambda}\frac{\arm(i,j)^2}{2h(i,j)}
&=\frac{n^{(p+1)/p}}{2}\iint_{D_{\gamma_n}}\frac{a_{\gamma_n}(x,y)^2}{a_{\gamma_n}(x,y)+n^{1/q-1/p}\ell_{\gamma_n}(x,y)+n^{-1/p}}\d x\d y.
\end{align*}
By use of the Lemmata~\ref{Lemma:symmdiff} and~\ref{Lemma:arm} it follows that
\begin{align*}
\abs{\iint_{D_{\gamma_n}}\frac{a_{\gamma_n}(x,y)^2}{a_{\gamma_n}(x,y)+n^{1/q-1/p}\ell_{\gamma_n}(x,y)+n^{-1/p}}\d x\d y
-\iint_{D_{\gamma}}a_{\gamma}(x,y)\d x\d y}
\to0
\end{align*}
as $n\to\infty$.
This establishes the asymptotic lower bound
\begin{align*}
C(\lambda^{(n)})
&>\frac{n^{(p+1)/p}}{2}\iint_{D_{\gamma}}a_{\gamma}(x,y)\d x\d y+o(n^{(p+1)/p})
&&\text{as }n\to\infty.
\end{align*}
In order to obtain an upper bound for the average case complexity, recall the algorithm for constructing the hook tableau during the application of the NPS algorithm described in \refs{algorithm}.
The reason why $\abs{H}$ can be less than $C(\lambda^{(n)})$ is that there might be a cancellation in step \textbf{H1}.
However, the total cancelation cannot exceed $\sum_{(i,j)\in\lambda^{(n)}}\leg(i,j)$ such that
\begin{align*}
C(\lambda^{(n)})
&<\E(\abs{H})+\sum_{(i,j)\in\lambda^{(n)}}\leg(i,j).
\end{align*}
Since the term
\begin{align*}
\sum_{(i,j)\in\lambda^{(n)}}\leg(i,j)
&=n^{(q+1)/q}\iint_{D_{\gamma_n}}\ell_{\gamma_n}(x,y)\d x\d y+\frac{n}{2}
\end{align*}
is of order less than $n^{(p+1)/p}$ as $n\to\infty$, we conclude \refq{Cp}.

Our starting point for the analysis of the worst case complexity is the inequality
\begin{align}\label{eq:Wineq}
\sum_{(i,j)\in\lambda}\arm(i,j)
&\leq W(\lambda^{(n)})
\leq\sum_{(i,j)\in\lambda^{(n)}}\arm(i,j)+\leg(i,j),
\end{align}
which is a trivial consequence of \refp{W}.
The right hand side of \refq{Wineq} equals
\begin{align*}
n^{(p+1)/p}\iint_{D_{\gamma_n}}a_{\gamma_n}(x,y)\d x\d y
+n^{(q+1)/q}\iint_{D_{\gamma_n}}\ell_{\gamma_n}(x,y)\d x\d y
+n.
\end{align*}
Again the term corresponding to the leg function is of lower order and can be dropped.
Thus $W(\lambda^{(n)})$ is asymptotically equivalent to the left hand side of \refq{Wineq}.
Lemmata~\ref{Lemma:symmdiff} and~\ref{Lemma:arm} imply that
\begin{align*}
\abs{\iint_{D_{\gamma_n}}a_{\gamma_n}(x,y)\d x\d y
-\iint_{D_{\gamma}}a_{\gamma}(x,y)\d x\d y}
\to0&&\text{as }n\to\infty,
\end{align*}
which yields \refq{Wp}.
The alternative formula in terms of hook coordinates is simply obtained by substitution, which completes the first part of the proof.
The second part, that is, the case $p>q$, follows similarly.
\end{proof}

\section{Partitions with two parts}\label{Section:zweizeiler}

The main result of this section is a nice formula for the average case complexity of the NPS algorithm when the partition consists of only two rows.

\begin{mythm}{zweizeiler} Let $\lambda=(\lambda_1,\lambda_2)$ be a partition with two parts. Then the average case complexity of the NPS algorithm on $\lambda$ is given by
\begin{align}\label{eq:zweizeiler}
C(\lambda)
&=\frac{\lambda_1(\lambda_1-1)}{4}+\frac{\lambda_2(\lambda_2-3)}{4}-2\sum_{k=1}^{\lambda_2}\binom{\lambda_2}{k}\frac{(-1)^k(2k-2)!}{(\lambda_1-\lambda_2+2)_{2k-1}}.
\end{align}
\end{mythm}

Our starting point for proving \reft{zweizeiler} is the following formula for the average case complexity of the NPS algorithm for general partitions.

\begin{mythm}{chicago}\textnormal{(\cite[Thm~2.4]{Sul2014})} Let $n\in\N$ and $\lambda$ be a partition of $n$. Then
\begin{align}\label{eq:chicago}
C(\lambda)
&=\sum_{x\in\lambda}\sum_{k=1}^{n}\abs{x}\frac{f^{\lambda}(x,k)}{f^{\lambda}}(H_n-H_{n-k}-1),
\end{align}
where the outer sum is taken over all cells $x$ of the Young diagram of $\lambda$, $\abs{x}=i+j-2$ denotes the distance of a cell $x=(i,j)$ to the top left cell, $f^{\lambda}$ denotes the number of SYT of shape $\lambda$, $f^{\lambda}(x,k)$ denotes the number of SYT of shape $\lambda$ such that the cell $x$ contains the entry $k$, and $H_n=\sum_{\ell=1}^{n}\frac{1}{\ell}$ denotes the $n$-th harmonic number.
\end{mythm}

While the appearance of harmonic numbers in \eqref{eq:chicago} is quite surprising, the most challenging expressions are the numbers $f^{\lambda}(x,k)$.
For partitions with only two parts we derive a first explicit expression for the average case complexity in terms of five double sums, each of which contains a harmonic number in the summand.

\begin{mylem}{sixdouble} Let $\lambda=(\lambda_1,\lambda_2)$ be a partition with two parts. Then
\begin{align*}
f^{\lambda}&=\frac{\lambda_1-\lambda_2+1}{\lambda_1+1}\binom{\lambda_1+\lambda_2}{\lambda_2}
\end{align*}
and
\begin{align}\label{eq:sixdouble}
C(\lambda)
&=\left(\binom{\lambda_1}{2}+\binom{\lambda_2+1}{2}\right)(H_{\lambda_1+\lambda_2}-1)\\
\notag&\quad-\frac{1}{f^{\lambda}}\sum_{j=1}^{\lambda_2}\sum_{k=j}^{2j-1}\frac{(j-1)(2j-k)}{k}\binom{k}{j}\binom{\lambda_1+\lambda_2-k}{\lambda_1-j}H_{\lambda_1+\lambda_2-k}\\
\notag&\quad+\frac{1}{f^{\lambda}}\sum_{j=1}^{\lambda_2}\sum_{k=j}^{2j-1}\frac{(j-1)(2j-k)}{k}\binom{k}{j}\binom{\lambda_1+\lambda_2-k}{\lambda_2-j-1}H_{\lambda_1+\lambda_2-k}\\
\notag&\quad-\frac{1}{f^{\lambda}}\sum_{j=\lambda_2+1}^{\lambda_1}\sum_{k=j}^{\lambda_2+j}\frac{(j-1)(2j-k)}{k}\binom{k}{j}\binom{\lambda_1+\lambda_2-k}{\lambda_1-j}H_{\lambda_1+\lambda_2-k}\\
\notag&\quad-\frac{1}{f^{\lambda}}\sum_{j=1}^{\lambda_2}\sum_{k=2j}^{\lambda_1+j}\frac{j(k-2j+2)}{k}\binom{k}{j-1}\binom{\lambda_1+\lambda_2-k}{\lambda_2-j}H_{\lambda_1+\lambda_2-k}\\
\notag&\quad+\frac{1}{f^{\lambda}}\sum_{j=1}^{\lambda_2}\sum_{k=2j}^{\lambda_2+j}\frac{j(k-2j+2)}{k}\binom{k}{j-1}\binom{\lambda_1+\lambda_2-k}{\lambda_1-j+1}H_{\lambda_1+\lambda_2-k}.
\end{align}
\end{mylem}

\begin{proof} The formula for $f^{\lambda}$ is a simple consequence of the hook-length formula~\refq{hlf}
Next note that
\begin{align*}
\sum_{k=1}^{n}\frac{f^{\lambda}(x,k)}{f^{\lambda}}=1,
\end{align*}
and therefore
\begin{align*}
\sum_{x\in\lambda}\sum_{k=1}^{n}\abs{x}\frac{f^{\lambda}(x,k)}{f^{\lambda}}(H_n-1)
&=(H_n-1)\sum_{x\in\lambda}\abs{x}
=(H_n-1)\left(\binom{\lambda_1}{2}+\binom{\lambda_2+1}{2}\right).
\end{align*}
Thus only the terms involving a harmonic numbers that depends on $k$ remain.
For a cell $x\in\lambda$ let $S(x)$ denote the set of all values $k\in[n]$ such that $f^{\lambda}(x,k)>0$.
We observe three cases,
\begin{align*}
S((1,j))&=\{j,\dots,2j-1\}&&\text{for }j\in[\lambda_2],\\
S((1,j))&=\{j,\dots,\lambda_2+j\}&&\text{for }j\in[\lambda_1]-[\lambda_2]\text{ and}\\
S((2,j))&=\{2j,\dots,\lambda_1+j\}&&\text{for }j\in[\lambda_2].
\end{align*}
Suppose $x=(1,j)$ is a cell of $\lambda$, $k\in S(x)$, and let $T$ be a SYT of shape $\lambda$ with $T(x)=k$.
Then the cells $y\in\lambda$ with $T(y)<k$ constitute the partition $\tilde{\mu}=(j-1,k-j)$.
On the other hand the cells $y\in\lambda$ with $T(y)>k$ form the skew shape $\lambda/\mu$ where $\mu=(j,k-j)$.
Moreover note that the set of SYT of shape $\lambda$ in which the cell $x$ contains the entry $k$ is in bijection with pairs of a SYT of shape $\tilde{\mu}$ and a SYT of skew shape $\lambda/\mu$.
Consequently
\begin{align*}
f^{\lambda}(x,k)
&=f^{\tilde{\mu}}f^{\lambda/\mu}.
\end{align*}
While $f^{\tilde{\mu}}$ is given by the hook-length formula, the number of SYT of skew shape $f^{\lambda/\mu}$ can be computed using Aitken's determinant formula~\cite{Aitken1943}
\begin{align*}
f^{\lambda/\mu}
&=(\abs{\lambda}-\abs{\mu})!\cdot\det_{i,j}\Big(1/(\lambda_i-\mu_j-i+j)!\Big).
\end{align*}
In the present case
\begin{align*}
f^{\lambda/\mu}
&=(\lambda_1+\lambda_2-k)!\cdot\det\begin{pmatrix}
1/(\lambda_1-j)! & 1/(\lambda_1+j-k+1)!\\
1/(\lambda_2-j-1)! & 1/(\lambda_2+j-k)!
\end{pmatrix}\\
&=\binom{\lambda_1+\lambda_2-k}{\lambda_1-j}-\binom{\lambda_1+\lambda_2-k}{\lambda_2-j-1}.
\end{align*}
Summing
\begin{align*}
\sum_{k\in S(x)}\abs{x}\frac{f^{\tilde{\mu}}f^{\lambda/\mu}(x,k)}{f^{\lambda}}(-H_{\lambda_1+\lambda_2-k})
\end{align*}
over all cells $x=(1,j)$ of the first row of $\lambda$, accounts for the first three double sums in \eqref{eq:sixdouble}.
The case where $x=(2,j)$ is a cell of the second row of $\lambda$ is treated in the same way, except that now $\tilde{\mu}=(k-j,j-1)$ and $\mu=(k-j,j)$, and accounts for the fourth and fifth double sums in \eqref{eq:sixdouble}.
\end{proof}

While the proof of \refl{sixdouble} is not too complicated, we begin to appreciate how remarkably simple the expression in~\eqref{eq:zweizeiler} really is, consisting of a single sum devoid of harmonic numbers. 

\medskip

In the following we will prove \reft{zweizeiler}. More precisely, denoting 
the right hand sides of~\eqref{eq:zweizeiler} and~\eqref{eq:sixdouble} by $A(\lambda_1,\lambda_2)$ and  $B(\lambda_1,\lambda_2)$, respectively, we
will show that 
\begin{equation}\label{Equ:A=B}
A(\lambda_1,\lambda_2)=B(\lambda_1,\lambda_2)
\end{equation}
holds for all $\lambda_1,\lambda_2\in\N$ with $0\leq\lambda_2\leq\lambda_1$.

Looking at the given problem, one could be tempted to try the following rather general summation tactic: compute for each of the sums a homogeneous recurrence relation in one of the discrete parameters, say $\lambda_2$ (using, e.g., the package \texttt{MultiSum}~\cite{Wegschaider:97}), and combine the found recurrences to one linear homogeneous recurrence for the expression 
\begin{equation}\label{Equ:ExprAMinusB}
T(\lambda_1,\lambda_2):=A(\lambda_1,\lambda_2)-B(\lambda_1,\lambda_2)
\end{equation}
(using, e.g., the Mathematica package~\texttt{GeneratingFunctions}~\cite{Mallinger:96}). However, in this particular situation this tactic seems rather clumsy: already the calculation of the linear recurrences for each single sum is a hard nut, and assembling the recurrences to a big recurrence for~\eqref{Equ:ExprAMinusB} is rather hopeless. 

Therefore we will follow an alternative tactic: Given such an expression in terms of complicated multi-sums, try to simplify the involved sums, and try to show that this identity holds in terms of these simpler objects. For definite hypergeometric sums this strategy has been worked out in the well known summation book $A=B$~\cite{AequalB}. E.g., suppose that we are given the definite sum on the left hand side of
\begin{equation}\label{Equ:HypId}
\sum_{i=1}^{\lambda_2} 2^{-i} \binom{i
+\lambda_2
}{i}=2^{\lambda_2}-1,\\
\end{equation}
which we denote by $S(\lambda_2)$. Then one can find the right hand side by computing a linear recurrence for $S(\lambda_2)$ using Zeilberger's summation paradigm of creative telescoping. We postpone any details and just present the found recurrence
$$S(\lambda_2+1)-2\,S(\lambda_2)=1.$$
Now one can read off the solution $2^{\lambda_2}-1$, and with the first initial value $S(1)=1$ the identity~\eqref{Equ:HypId} is established. For more involved examples, one might start with a hypergeometric sum and obtains a linear recurrence of higher order (and not just of order one as above). Then one can use in addition, e.g., Petkov{\v s}ek's recurrence solver~\cite{AequalB} to find all hypergeometric solutions. In the end, one can possibly express the definite sum in terms of these objects. 

\medskip

However, the types of sums arising in~\eqref{Equ:A=B} are not as simple as the one in~\eqref{Equ:HypId}, in particular the hypergeometric technology from~\cite{AequalB} is not general enough to prove identity~\eqref{Equ:A=B}. To overcome this situation, we will utilize the summation algorithms in the setting of difference rings~\cite{Karr:81,DR1,DR2}. Namely, instead of working only with hypergeometric products, we will apply symbolic summation tools that are tuned for expressions in terms indefinite nested sums defined over hypergeometric products. This more general class of summation objects can be defined as follows.

\medskip

\noindent\textbf{Definition.} Let $f(n)$ be an expression that evaluates at non-negative 
integers (from a certain point on) to elements of a field $\K$ containing the rational numbers $\Q$. Then $f(n)$ 
is called an \emph{expression in terms of indefinite nested sums over hypergeometric products  w.r.t.\ $n$}  if it is composed by elements from the rational function field $\K(n)$, by the three operations 
($+,-,\cdot$), by \emph{hypergeometric products} of the form $\prod_{k=l}^nh(k)$ with $l\in\N$ and a rational function $h(k)\in\K(k)\setminus\{0\}$, and by sums of the form $\sum_{k=l}^nF(k)$ with $l\in\N$ and where 
$F(k)$, being free of $n$, is an expression in terms of indefinite nested sums over hypergeometric products w.r.t.\ $k$. 

\medskip

More precisely, the summation package \texttt{Sigma}~\cite{Sigma} based on the algorithmic difference ring theory~\cite{Karr:81,DR1,DR2} can tackle the following definite summation problem.

\medskip

\begin{center}
\fbox{\begin{minipage}{15.5cm}
\textbf{\textbf{Problem T}: Transformation of a definite sum to indefinite nested sums.}\\ 
\small
\textit{Given} a definite sum, say $S(n)=\sum_{k=0}^nf(n,k)$, where $f(n,k)$ is given in terms of indefinite nested sums over hypergeometric products w.r.t.\footnote{Here $n$ is considered as a parameter which does not occur in any summation bound of the indefinite nested sums of $f(n,k)$.}\ $k$\\
\textit{find} an expression $T(n)$ in terms of indefinite nested sums over hypergeometric products w.r.t.\ $n$ and \textit{find} a $\lambda\in\N$ such that $S(\nu)=T(\nu)$ holds for all $\nu\in\N$ with $\nu\geq\lambda$.
\end{minipage}}
\end{center}
\medskip

\noindent In order to tackle Problem~T, the following summation steps can be carried out in \texttt{Sigma}.
\begin{enumerate}
 \item Compute a linear recurrence of $S(n)$ in $n$, say of order $d$.
 \item Solve the recurrence in terms of d'Alembertian solutions~\cite{Abramov:94}, i.e., in terms of all solutions that are expressible in terms of indefinite nested sums over hypergeometric products w.r.t.\ $n$.
 \item Combine the derived solutions yielding an expression  $T(n)$ in terms of indefinite nested sums over hypergeometric products such that $S(\nu)=T(\nu)$ holds for $\nu=\lambda,\lambda+1,\dots \lambda+d$ for some appropriately chosen $\lambda\in\N$. 
\end{enumerate}
\noindent If one succeeds in this strategy, it follows with some mild side conditions that $S(\nu)=T(\nu)$ holds for all $\nu\geq\lambda$. In other words, one has solved Problem~T. Otherwise, if one succeeds in computing a recurrence, but fails to combine the solutions accordingly, it follows that there does not exist such a representation of $S(n)$ in terms of indefinite nested sums.

\medskip

We remark that not any definite sum as specified in Problem~T can be transformed to a representation in terms of indefinite nested sums over hypergeometric products. But as it will turn out below, \texttt{Sigma}'s summation toolbox, in particular its solution for Problem~T, can be used iteratively to handle the multi-sums arising in~\eqref{eq:zweizeiler}.

Now suppose that we derived an alternative representation of~\eqref{Equ:ExprAMinusB} in terms of our class of indefinite nested sums defined over hypergeometric products. Then we will be in the position to utilize the following very special feature~\cite[Prop.~7.3]{DR2}.

\medskip

\begin{center}
\fbox{\begin{minipage}{15.5cm}
\textbf{\textbf{Problem S}: Simplification of indefinite nested sums.}\\ 
\small
\textit{Given} an expression $T(n)$ in terms of indefinite nested sums over hypergeometric products;\\ 
\textit{find} an expression $\tilde{T}(n)$ in terms of indefinite nested sums over hypergeometric products and \textit{find} a $\delta\in\N$ with the following properties:
\begin{enumerate}
 \item $T(\nu)=\tilde{T}(\nu)$ for all $\nu\in\N$ with $\nu\geq\delta$;
 \item the nested sums and hypergeometric products in $\tilde{T}(k)$ (except products of the form $\alpha^k$ with $\alpha$ being a root of unity) are algebraically independent among each other.
 \end{enumerate}
 \end{minipage}}
\end{center}

 \medskip
 
\noindent We emphasize that such a computed $\tilde{T}(n)$ has the following special property: If $T(\nu)=0$ holds for all $\nu\geq\delta$ for some $\delta\in\N$, then $\tilde{T}(n)$ is the zero-expression (or it can be simplified to $0$ by simple polynomial arithmetic). Precisely this will happen to our expression in~\eqref{Equ:ExprAMinusB}, after we transformed it to indefinite nested sums and eliminated all algebraic relations among the arising summation objects. 

\medskip

Summarizing, we will prove that $T(\lambda_1,\lambda_2)$ from~\eqref{Equ:ExprAMinusB} evaluates to zero for any $\lambda_1,\lambda_2\in\N$ with $0\leq\lambda_2\leq\lambda_1$ by executing the following two main steps.

\medskip

\begin{description}
\item[(DEF)] Using \texttt{Sigma's} definite summation toolbox (see Problem~T), we will find alternative sum representations where the occurring sums are indefinite nested w.r.t.\ to the discrete parameter $\lambda_2$. In a nutshell, we will rewrite the expression $T(\lambda_1,\lambda_2)$ given in terms of $6$ definite sums to an expression in terms of indefinite nested sums w.r.t.\ $\lambda_2$.

\item[(IND)] Using \texttt{Sigma's} indefinite summation toolbox (see Problem~S), we will rewrite the expression~\eqref{Equ:ExprAMinusB} further such that no algebraic relations exist among the arising indefinite nested sums and products. As we will see below, the derived expression of $\tilde{T}(\lambda_1,\lambda_2)$ will collapse to zero, which will prove~\eqref{Equ:A=B} and thus will establish \reft{zweizeiler}.

\end{description}

\medskip

We start our proposed summation tactic by loading in the package

\begin{mma}
\In << Sigma.m \\
\Print \LoadP{Sigma - A summation package by Carsten Schneider
\copyright\ RISC-Linz}\\
\end{mma}

\noindent into the Mathematica system. First we tackle the definite sum 
$$S_0(\lambda_1,\lambda_2)=\sum_{k=1}^{\lambda_2}\binom{\lambda_2}{k}\frac{(-1)^k(2k-2)!}{(\lambda_1-\lambda_2+2)_{2k-1}}$$
given in $A(\lambda_1,\lambda_2)$, i.e., given in the right hand side of~\eqref{eq:zweizeiler}. After entering this sum into Mathematica
\begin{mma}
 \In S0=SigmaSum[\frac{(-1)^k SigmaBinomial[\lambda_2,k] (-2+2 k)!}{SigmaPochhammer[
        2
        +\lambda_1
        -\lambda_2,{-1+2 k}]},\{k,1,\lambda_2\}]\\
\Out \sum_{k=1}^{\lambda_2} \frac{(-1)^k \binom{\lambda_2}{k} (-2+2 k)!}{\big(
        2
        +\lambda_1
        -\lambda_2
\big)_{-1+2 k}}\\
\end{mma}

\noindent we use \texttt{Sigma}'s recurrence finder to calculate the following recurrence relation:

\begin{mma}\MLabel{MMA:GR}
\In rec = GenerateRecurrence[S0,\lambda_2][[1]]\\ 
\Out 2 \big(
        \lambda_2+1\big) \text{SUM}[\lambda_2]
+\big(
        -3
        -\lambda_1
        -3 \lambda_2
\big) \text{SUM}[\lambda_2+1]
+\big(
        2
        +\lambda_1
        +\lambda_2
\big) \text{SUM}[\lambda_2+2]
=\frac{-\lambda_1
-\lambda_2
-2
}{\lambda_1
-\lambda_2
}\\  
\end{mma}
\noindent This means that $S_0(\lambda_1,\lambda_2)=\texttt{S0}=\texttt{SUM}[\lambda_2]$ is a solution of~\myOut{\ref{MMA:GR}}. Internally, the creative telescoping paradigm is used: given the summand
$$f(\lambda_1,\lambda_2,k)=\binom{\lambda_2}{k}\frac{(-1)^k(2k-2)!}{(\lambda_1-\lambda_2+2)_{2k-1}}$$
of $S_0(\lambda_1,\lambda_2)$,
\texttt{Sigma} computes the constants 
$c_0(\lambda_1,\lambda_2)=2 \big(
                \lambda_2+1\big)$, $c_1(\lambda_1,\lambda_2)=-3
        -\lambda_1
        -3 \lambda_2$ and 
        $c_2(\lambda_1,\lambda_2)=2
        +\lambda_1
        +\lambda_2$ together with the expression
$$ g(\lambda_1,\lambda_2,k)=\tfrac{k \big(\lambda_2+1\big)\big(
                2 k
                +\lambda_1
                -\lambda_2
        \big)
\big(-1
                +2 k
                +\lambda_1
                -\lambda_2
        \big)
\big(2
                +\lambda_1
                +\lambda_2
        \big)}{\big(
                -2
                +k
                -\lambda_2
        \big)
\big(-1
                +k
                -\lambda_2
        \big)
\big(-\lambda_1
                +\lambda_2
        \big)
\big(-1
                -\lambda_1
                +\lambda_2
        \big)
} 
\frac{(-1)^{k-1} \binom{\lambda_2}{k} (2 k-2)!}{\big(2
                +\lambda_1
                -\lambda_2
        \big)_{2 k-1}}
$$
such that the summand recurrence
\begin{multline}\label{Equ:SummandRec}
c_0(\lambda_1,\lambda_2)f(\lambda_1,\lambda_2,k)+c_1(\lambda_1,\lambda_2)f(\lambda_1,\lambda_2+1,k)+c_2(\lambda_1,\lambda_2)f(\lambda_1,\lambda_2+2,k)\\
=g(\lambda_1,\lambda_2,k+1)-g(\lambda_1,\lambda_2,k)
\end{multline}
holds. Note that the special case $\lambda_1=\lambda_2$ is problematic and will be treated separately. For all other cases, i.e., for all $\lambda_1,\lambda_2,k\in\N$ with $\lambda_1>\lambda_2\geq k\geq0$ the correctness of this relation can be checked easily by simple polynomial arithmetic. Hence summing~\eqref{Equ:SummandRec} over $k$ yields the recurrence~\myOut{\ref{MMA:GR}}. In particular, since~\eqref{Equ:SummandRec} has been verified, we proved that $S_0(\lambda_1,\lambda_2)$ is a solution of~\myOut{\ref{MMA:GR}} for all $\lambda_1,\lambda_2\in\N$ with $\lambda_1>\lambda_2$.
We remark further that in this particular instance also simpler algorithms, like the Mathematica implementation~\cite{Paule:95} of Zeilberger's creative telescoping algorithm~\cite{AequalB} for hypergeometric terms, could have been used.

Next, we exploit \texttt{Sigma}'s recurrence solver to compute all d'Alembertian solutions~\cite{Abramov:94,AequalB} by executing the function call

\begin{mma}
\In recSol = SolveRecurrence[rec, SUM[\lambda_2]] // ToSpecialFunction\\
\Out\{
\{0 , -1
        -\lambda_1
        +\lambda_2\},
         \{0 , \frac{2^{\lambda_2} \lambda_2!}{\big(
                \lambda_1+2\big)_{\lambda_2}}
        +\frac{1}{2} \big(
                1
                +\lambda_1
                -\lambda_2
        \big) 
        \sum_{i=1}^{\lambda_2} \frac{2^i i!}{i \big(
                2+\lambda_1\big)_i}\},
        \{ 1 , -\big(
                -1
                -\lambda_1
                +\lambda_2
        \big) 
        \sum_{i=1}^{\lambda_2} \frac{1}{1
        -i
        +\lambda_1
        }\newline
        +\frac{1}{2} \big(
                -1
                -\lambda_1
                +\lambda_2
        \big) 
        \sum_{i=1}^{\lambda_2} \frac{2^i i! 
        \sum_{j=1}^i \frac{2^{-j} \big(
                2+\lambda_1\big)_j}{j!}}{i \big(
                2+\lambda_1\big)_i}
        -\frac{2^{\lambda_2} \lambda_2!}{\big(
                \lambda_1+2\big)_{\lambda_2}} 
        \sum_{i=1}^{\lambda_2} \frac{2^{-i} \big(
                2+\lambda_1\big)_i}{i!}
        +\frac{1}{2} \big(
                1
                +\lambda_1
                -\lambda_2
        \big) 
        \sum_{i=1}^{\lambda_2} \frac{1}{i}\}\}\\
\end{mma}
\noindent This means that  
\begin{align*}
h_1(\lambda_1,\lambda_2)&=-1
        -\lambda_1
        +\lambda_2,\\
h_2(\lambda_1,\lambda_2)&=\frac{2^{\lambda_2} \lambda_2!}{\big(
                \lambda_1+2\big)_{\lambda_2}}
        +\frac{1}{2} \big(
                1
                +\lambda_1
                -\lambda_2
        \big) 
        \sum_{i=1}^{\lambda_2} \frac{2^i i!}{i \big(
                2+\lambda_1\big)_i}
\end{align*} 
are two linearly independent solutions of the homogeneous version of the recurrence~\myOut{\ref{MMA:GR}} and that
\begin{multline*}
p(\lambda_1,\lambda_2)= \big(
                -1
                -\lambda_1
                +\lambda_2
        \big) \Big(
        -\sum_{i=1}^{\lambda_2} \frac{1}{1
        -i
        +\lambda_1
        }-\sum_{i=1}^{\lambda_2} \frac{1}{i}
        +\frac{1}{2} 
        \sum_{i=1}^{\lambda_2} \frac{2^i i! 
        \sum_{j=1}^i \frac{2^{-j} \big(
                2+\lambda_1\big)_j}{j!}}{i \big(
                2+\lambda_1\big)_i}\Big)\\
        -\frac{2^{\lambda_2} \lambda_2!}{\big(
                \lambda_1+2\big)_{\lambda_2}} 
        \sum_{i=1}^{\lambda_2} \frac{2^{-i} \big(
                2+\lambda_1\big)_i}{i!}
\end{multline*}
is a particular solution of the recurrence~\myOut{\ref{MMA:GR}} itself. 
In other words, letting $d_1(\lambda_1)$ and $d_2(\lambda_1)$ be arbitrary constants we obtain the general solution 
\begin{equation}\label{Equ:GeneralSol}
d_1(\lambda_1)\,h_1(\lambda_1,\lambda_2)+
        d_2(\lambda_1)\,h_2(\lambda_1,\lambda_2)+p(\lambda_1,\lambda_2)
\end{equation}
of the recurrence~\myOut{\ref{MMA:GR}}. By the underlying algorithms the d'Alembertian solutions are produced in a rather complicated form and \texttt{Sigma}'s built in simplifier worked hard to drop the nice solutions.
In addition, the arising sums within the output are algebraically independent among each other. For further details on these aspects we refer to~\cite{Schneider:13a,DR1,DR2} and references therein. 
We emphasize further that one can verify straightforwardly with simple polynomial arithmetic that~\eqref{Equ:GeneralSol} is indeed a solution of the recurrence~\myOut{\ref{MMA:GR}} for all $\lambda_1,\lambda_2\in\N$ with $\lambda_1>\lambda_2$. 

Finally, we determine the values
$d_1(\lambda_1)=-\frac{1}{\lambda_1+1}$ and $d_2(\lambda_1)=-1$ such that $S_0(\lambda_1,\lambda_2)$ and~\eqref{Equ:GeneralSol} agree for $\lambda_2=1,2$. Since both expressions are a solution of~\myOut{\ref{MMA:GR}} by construction, it follows that they agree for all $\lambda_2\in\N$. This last step can be performed by computing the first two initial values at $\lambda_2=1,2$, namely
\begin{mma}
 \In initial=Table[S0,\{\lambda_2,1,2\}]\\
 \Out \big\{-\frac{1}{\lambda_1+1},-\frac{2 \big(
        \lambda_1^2+3 \lambda_1+1\big)}{\lambda_1 \big(
        \lambda_1+1
\big)
\big(\lambda_1+2\big)}\big\}\\
 \end{mma}
\noindent and  using \texttt{Sigma}'s function call to combine the solutions accordingly:
\begin{mma}\MLabel{MMA:FLC}
 \In FindLinearCombination[recSol, \{1, initial\}, \lambda_2, 2]\\
 \Out \frac{1}{2} \big(
        -1
        -\lambda_1
        +\lambda_2
\big) 
\sum_{i=1}^{\lambda_2} \frac{2^i i! 
\sum_{j=1}^i \frac{2^{-j} \big(
        2+\lambda_1\big)_j}{j!}}{i \big(
        2+\lambda_1\big)_i}
+\frac{1}{2} \big(
        -1
        -\lambda_1
        +\lambda_2
\big) 
\sum_{i=1}^{\lambda_2} \frac{2^i i!}{i \big(
        2+\lambda_1\big)_i}
-\frac{2^{\lambda_2} \lambda_2!}{\big(
        \lambda_1+2\big)_{\lambda_2}} 
\sum_{i=1}^{\lambda_2} \frac{2^{-i} \big(
        2+\lambda_1\big)_i}{i!}
-\frac{2^{\lambda_2} \lambda_2!}{\big(
        \lambda_1+2\big)_{\lambda_2}}
+\big(
        1
        +\lambda_1
        -\lambda_2
\big) 
\sum_{i=1}^{\lambda_2} \frac{1}{1
-i
+\lambda_1
}
+\frac{1}{2} \big(
        1
        +\lambda_1
        -\lambda_2
\big) 
\sum_{i=1}^{\lambda_2} \frac{1}{i}
+\frac{1
+\lambda_1
-\lambda_2
}{\lambda_1+1}
\\
\end{mma}
\noindent Summarizing, we showed that $S_0(\lambda_1,\lambda_2)$ equals the expression calculated in~\myOut{\ref{MMA:FLC}} for all $\lambda_1,\lambda_2\in\N$ with $1\leq\lambda_2<\lambda_1$.

\medskip

This ``hand calculation'' with the computer is reasonable if one treats the simplest of the sums in~\eqref{Equ:ExprAMinusB}. However, if one is faced with one of the double sums in~\eqref{eq:sixdouble}, this mechanical task gets more and more tedious. Luckily, we can use the package \texttt{EvaluateMultiSums}~\cite{Schneider:13a} that has been originally designed for similar and even worse expressions arising in the field of particle physics, see~\cite{Schneider:16b} and references therein. Namely, after loading in the package
\begin{mma}
\In << EvaluateMultiSums.m \\
\Print \LoadP{EvaluateMultiSums by Carsten Schneider
\copyright\ RISC-Linz}\\
\end{mma}

\noindent we can carry out the above calculations completely automatically using internally the functionality of \texttt{Sigma}. For instance, typing in the command
\begin{mma}\MLabel{MMA:solSum0Pre}
\In solSum0=EvaluateMultiSum[\frac{(-1)^k \binom{\lambda_2}{k} (-2+2 k)!}{\big(
        2
        +\lambda_1
        -\lambda_2
\big)_{-1+2 k}},\{\{k,1,\lambda_2\}\},\{\lambda_2,\lambda_1\},\{1,1\},\{\lambda_1,\infty\}]\\ 
\Out  \frac{1}{2} \big(
        -1
        -\lambda_1
        +\lambda_2
\big) 
\sum_{i=1}^{\lambda_2} \frac{2^i i! 
\sum_{j=1}^i \frac{2^{-j} \big(
        2+\lambda_1\big)_j}{j!}}{i \big(
        2+\lambda_1\big)_i}
+\frac{1}{2} \big(
        -1
        -\lambda_1
        +\lambda_2
\big) 
\sum_{i=1}^{\lambda_2} \frac{2^i i!}{i \big(
        2+\lambda_1\big)_i}
-\frac{2^{\lambda_2} \lambda_2!}{\big(
        \lambda_1+2\big)_{\lambda_2}} 
\sum_{i=1}^{\lambda_2} \frac{2^{-i} \big(
        2+\lambda_1\big)_i}{i!}
-\frac{2^{\lambda_2} \lambda_2!}{\big(
        \lambda_1+2\big)_{\lambda_2}}
+\big(
        1
        +\lambda_1
        -\lambda_2
\big) 
\sum_{i=1}^{\lambda_2} \frac{1}{1
-i
+\lambda_1
}
+\frac{1}{2} \big(
        1
        +\lambda_1
        -\lambda_2
\big) 
\sum_{i=1}^{\lambda_2} \frac{1}{i}
+\frac{1
+\lambda_1
-\lambda_2
}{\lambda_1+1}
\\
\end{mma}
\noindent we arrive at an equivalent result \texttt{solSum0} in terms of factorials and the Pochhammer symbol. Within the function call the input $\{\lambda_2,\lambda_1\},\{1,1\},\{\lambda_1,\infty\}$ specifies that the discrete parameters
$\lambda_1,\lambda_2$ fulfill the constraints $1\leq\lambda_2\leq\lambda_1$ and $1\leq\lambda_1\leq\infty$. 

\noindent Finally, we rewrite the objects in \texttt{solSum0} in terms of the harmonic numbers and the binomial coefficient. This rewriting can be accomplished by the function call 
\begin{mma}\MLabel{MMA:solSum0}
\In solSum0=SigmaReduce[solSum0,\lambda_2,Tower\to \big\{\binom{\lambda_1
        +\lambda_2
        }{\lambda_2},H_{\lambda_1
        +\lambda_2
        },H_{\lambda_1
        -\lambda_2
        },H_{\lambda_2}\big\}] \\
\Out \frac{2^{\lambda_2}}{\binom{\lambda_1
+\lambda_2
}{\lambda_2}} \big(
        -
        \sum_{i=1}^{\lambda_2} 2^{-i} \binom{i
        +\lambda_1
        }{i}-1\big)
+\frac{1}{2} \big(
        -1
        -\lambda_1
        +\lambda_2
\big) 
\sum_{i=1}^{\lambda_2} \frac{2^i 
\sum_{j=1}^i 2^{-j} \binom{j
+\lambda_1
}{j}}{i \binom{i
+\lambda_1
}{i}}
+\frac{1}{2} \big(
        -1
        -\lambda_1
        +\lambda_2
\big) 
\sum_{i=1}^{\lambda_2} \frac{2^i}{i \binom{i
+\lambda_1
}{i}}
+\big(
        -1
        -\lambda_1
        +\lambda_2
\big) H_{\lambda_1
-\lambda_2
}
+\big(
        1
        +\lambda_1
        -\lambda_2
\big) H_{\lambda_1}
+\frac{1}{2} \big(
        1
        +\lambda_1
        -\lambda_2
\big) H_{\lambda_2}
+1\\ 
\end{mma}

\begin{rem}\label{Remark:Rewriting} 
Internally, the underlying difference ring machinery of \texttt{Sigma} is activated~\cite{Schneider:13a}. Loosely speaking, the Pochhammer symbols and harmonic numbers arising in~\myOut{\ref{MMA:solSum0Pre}} are rewritten in terms of $\tbinom{\lambda_1+\lambda_2}{\lambda_2}$, and as a consequence also the arising sums are then rephrased in terms of $\tbinom{\lambda_1+i}{i}$ and $\tbinom{\lambda_1+j}{j}$. Here any polynomial in terms of the summation variables $i$ or $j$ that occurs in the summands of the numerators or denominators is expanded using partial fraction decomposition. This might lead to several new sums with simpler summands. Inside of this construction also Problem~S is carried out: a subset of these sums is taken whose elements are algebraically independent among each other, and the other sums are rewritten in terms of these algebraic independent sums. This finally leads to the output in~\myOut{\ref{MMA:solSum0}}.
\end{rem}

\noindent Summarizing, we end up at the following alternative representation
\begin{multline}\label{Equ:S0Rep}
S_0(\lambda_1,\lambda_2)=-\frac{2^{\lambda_2}}{\binom{\lambda_1
+\lambda_2
}{\lambda_2}} \Big(
        \sum_{i=1}^{\lambda_2} \frac{\binom{i
        +\lambda_1
        }{i}}{2^{i}}+1\Big)
+\frac{-1
        -\lambda_1
        +\lambda_2}{2}  
\sum_{i=1}^{\lambda_2} \frac{\displaystyle 2^i 
\sum_{j=1}^i\tfrac{\binom{j
+\lambda_1
}{j}}{2^{j}}}{i \binom{i
+\lambda_1
}{i}}\\
-\frac{1
        +\lambda_1
        -\lambda_2}{2}
\sum_{i=1}^{\lambda_2} \frac{2^i}{i \binom{i
+\lambda_1
}{i}}
+\big(  1
        +\lambda_1
        -\lambda_2
\big)\big( 
+H_{\lambda_1}
+\frac{1}{2}H_{\lambda_2}-H_{\lambda_1
-\lambda_2
}\big)
+1
\end{multline}
where all the sums are indefinite nested w.r.t.\ the outer most summation index $\lambda_2$.

\medskip

Finally, we have to address the special case $\lambda_1=\lambda_2$. In this particular case the sum $S_0(\lambda_2,\lambda_2)$ simplifies to

\begin{mma}\MLabel{MMA:solSum0Equal}
\In solSum0Equal=EvaluateMultiSum[\frac{(-1)^k \binom{\lambda_2}{k} (2 k-2)!}{(2)_{2 k-1}},\{\{k,1,\lambda_2\}\},\{\lambda_2\},\{1\},\{\infty\}]\\ 
\Out -\frac{2^{2 \lambda_2}\lambda_2!^2}{(2 \lambda_2)!}
+\frac{H_{\lambda_2}}{2}
+1\\
\end{mma}
\noindent which can be easily rewritten in terms of the binomial coefficient:
\begin{mma}
\In SigmaReduce[solSum0Equal,\lambda_2,Tower\to\{\tbinom{2\lambda_2}{\lambda_2}]\\
\Out -\frac{2^{2 \lambda_2}}{\binom{2 \lambda_2}{\lambda_2}}
+\frac{H_{\lambda_2}}{2}
+1\\
\end{mma}

\noindent Hence we calculated the nice simplification
\begin{equation}\label{Equ:S0SpecialCase}
S_0(\lambda_2,\lambda_2)=-\frac{2^{2 \lambda_2}}{\binom{2 \lambda_2}{\lambda_2}}
+\frac{H_{\lambda_2}}{2}
+1.
\end{equation}
Note further that our result \texttt{solSum0}, i.e., the right hand side of~\eqref{Equ:S0Rep} is a well defined expression for $\lambda_1=\lambda_2$. More precisely, we get

\begin{mma}
\In solSum0Subst=solSum0/.\lambda_1\to\lambda_2\\
\Out -\frac{1}{2} 
\sum_{i=1}^{\lambda_2} \frac{2^i 
\sum_{j=1}^i 2^{-j} \binom{j
+\lambda_2
}{j}}{i \binom{i
+\lambda_2
}{i}}
-\frac{2^{\lambda_2}}{\binom{2 \lambda_2}{\lambda_2}} 
\sum_{i=1}^{\lambda_2} 2^{-i} \binom{i
+\lambda_2
}{i}
-\frac{1}{2} 
\sum_{i=1}^{\lambda_2} \frac{2^i}{i \binom{i
+\lambda_2
}{i}}
-\frac{2^{\lambda_2}}{\binom{2 \lambda_2}{\lambda_2}}
+\frac{3 H_{\lambda_2}}{2}
+1\\
\end{mma}
\noindent Observe further that the arising sums in \texttt{solSum0Subst} can be simplified with our machinery (see Problem~T). Namely, for all $\lambda_2\in\N$ we calculate the right hand sides of the identities given in~\eqref{Equ:HypId} and 
\begin{align*}
\sum_{i=1}^{\lambda_2} \frac{2^i}{i \binom{i
+\lambda_2
}{i}}&=2^{-\lambda_2} 
\sum_{i=1}^{\lambda_2} \frac{2^i}{i},\\
\sum_{i=1}^{\lambda_2} \frac{2^i}{i \binom{i
+\lambda_2
}{i}}\sum_{j=1}^i 2^{-j} \binom{j
+\lambda_2
}{j}&=-2^{-\lambda_2} 
\sum_{i=1}^{\lambda_2} \frac{2^i}{i}
+2 H_{\lambda_2}.
\end{align*}
Inserting these simplifications in \texttt{solSum0Subst} we arrive at the same result as given in~\myOut{\ref{MMA:solSum0Equal}}. Note that the described calculations can be carried out straightforwardly by simply typing in
\begin{mma}
\In EvaluateMultiSum[solSum0Equal,\{\},\{\lambda_2\},\{1\},\{\infty\}];\\ 
\Out -\frac{2^{2 \lambda_2}}{\binom{2 \lambda_2}{\lambda_2}}
+\frac{H_{\lambda_2}}{2}
+1\\
\end{mma}
\noindent Since this result is precisely the same expression as given on the right hand side  of~\eqref{Equ:S0SpecialCase}, it follows that our identity~\eqref{Equ:S0Rep} holds also for the special case $\lambda_2=\lambda_1$, i.e., it holds for all $\lambda_1,\lambda_2\in\N$ with $1\leq\lambda_2\leq\lambda_1$.

\smallskip

Next, we turn to the first double sum 
$$S_1(\lambda_1,\lambda_2)=\sum_{j=1}^{\lambda_2}\sum_{k=j}^{2j-1}\frac{(j-1)(2j-k)}{k}\binom{k}{j}\binom{\lambda_1+\lambda_2-k}{\lambda_1-j}H_{\lambda_1+\lambda_2-k}
$$
of the expression $B(\lambda_1,\lambda_2)$, i.e., of the right hand side of~\eqref{eq:sixdouble}.
For convenience, we adapt the sum so that the lower bounds of the summation quantifiers are non-negative integers. This yields 
\begin{equation}\label{Equ:OuterSum}
S_1(\lambda_1,\lambda_2)=
\sum_{j=1}^{\lambda_2} h(\lambda_1,\lambda_2,j)
\end{equation}
with
$$h(\lambda_1,\lambda_2,j)=\sum_{k=0}^{-1+j} \frac{(-1+j)(j
        -k
        )}{j
        +k
        }\binom{j
        +k
        }{j} \binom{-j
        -k
        +\lambda_1
        +\lambda_2
        }{-j
        +\lambda_1
        }H_{-j
        -k
        +\lambda_1
        +\lambda_2
        }.$$
Here we apply iteratively our summation machinery for Problem~T.
We start with the inner sum $h(\lambda_1,\lambda_2,j)$. First, we compute a linear recurrence of $h(\lambda_1,\lambda_2,j)$ in $j$ of order 2, and afterwards we compute two linearly independent solutions plus one particular solution of the found recurrence. Finally, taking the first two initial values yields a rather huge expression (which is too big to be printed here) in terms of indefinite nested sums w.r.t.\ the integer parameter $j$. Precisely this form enables one to apply \texttt{Sigma}'s summation toolbox again to solve Problem~T for the definite sum~\eqref{Equ:OuterSum} where $j$ is the main summation variable: we can calculate a linear recurrence of~\eqref{Equ:OuterSum} in $\lambda_2$, solve the recurrence in terms of d'Alembertian solutions, and obtain finally a representation where the occurring sums are indefinite nested w.r.t.\ the parameter $\lambda_2$. In order to carry out all the calculations steps automatically, we execute simply the function call
\begin{mma}
 \In solSum1=EvaluateMultiSum[\frac{(-1+j) \binom{j
        +k
        }{j} \binom{-j
        -k
        +\lambda_1
        +\lambda_2
        }{-j
        +\lambda_1
        } (j
        -k
        ) H_{-j
        -k
        +\lambda_1
        +\lambda_2
        }}{j
        +k
        },\newline
        \hspace*{5cm}\{\{k,0,j\},\{j,1,\lambda_2\}\},\{\lambda_2,\lambda_1\},\{1,1\},\{\lambda_1,\infty\}];\\
\end{mma}
\noindent and get the result in terms of Pochhammer symbols and factorials. As explained in Remark~\ref{Remark:Rewriting} this result (which we did not print here) can be rewritten in terms of the binomial coefficient $\binom{\lambda_1+\lambda_2}{\lambda_2}$ by executing the following command:
\begin{mma}
\In solSum1=SigmaReduce[solSum1,\lambda_2,Tower\to\{\tbinom{\lambda_1+\lambda_2}{\lambda_2}\}]\\
 \Out 
 \binom{\lambda_1
+\lambda_2
}{\lambda_2} 
\frac{1}{16}\big(7 \lambda_1
        -3 \lambda_1^2
        -7 \lambda_2
        +4 \lambda_1 \lambda_2
        -\lambda_2^2
\big) \Big(
        -
        \sum_{i=1}^{\lambda_2} \frac{2^i 
        \sum_{j=1}^i 2^{-j} \binom{j
        +\lambda_1
        }{j}}{i \binom{i
        +\lambda_1
        }{i}}
        +H_{\lambda_2}
        -
        \sum_{i=1}^{\lambda_2} \frac{2^i}{i \binom{i
        +\lambda_1
        }{i}}
\Big)
\newline
+2^{\lambda_2-4}\big(-8
        +3 \lambda_1
        +\lambda_1^2
        -\lambda_2
        -2 \lambda_1 \lambda_2
        +\lambda_2^2
\big)\Bigg(
        H_{\lambda_1} \big(
                1-2 
                \sum_{i=1}^{\lambda_2} 2^{-i} \binom{i
                +\lambda_1
                }{i}\big)
\newline
\hspace*{0.4cm}  +\big(
                \sum_{i=1}^{\lambda_2} \frac{2^i}{i \binom{i
                +\lambda_1
                }{i}}\big) 
        \sum_{i=1}^{\lambda_2} 2^{-i} \binom{i
        +\lambda_1
        }{i}
   +\frac{2^{\lambda_2+1}}{\binom{\lambda_1
        +\lambda_2
        }{\lambda_2} \Big(
                1
                +\lambda_1
                +\lambda_2
        \Big)} \big(
                \sum_{i=1}^{\lambda_2} 2^{-i} \binom{i
                +\lambda_1
                }{i}\big)^2
        +2 
        \sum_{i=1}^{\lambda_2} 2^{-i} \binom{i
        +\lambda_1
        }{i} H_{-i
        +\lambda_1
        }
\newline
\hspace*{0.4cm}         -
        \sum_{i=1}^{\lambda_2} \frac{2^i 
        \sum_{j=1}^i 2^{-j} \binom{j
        +\lambda_1
        }{j}}{i \binom{i
        +\lambda_1
        }{i}}
      +\big(
                \sum_{i=1}^{\lambda_2} 2^{-i} \binom{i
                +\lambda_1
                }{i}\big) 
        \sum_{i=1}^{\lambda_2} \frac{2^i 
        \sum_{j=1}^i 2^{-j} \binom{j
        +\lambda_1
        }
        {j}}{i \binom{i
        +\lambda_1
        }{i}}
      -2 
        \sum_{i=1}^{\lambda_2} 
        \frac{2^i \big(
                \sum_{j=1}^i 2^{-j} \tbinom{j
                +\lambda_1
                }{j}\big)^2}{\tbinom{i
        +\lambda_1
        }{i} \big(
                1
                +i
                +\lambda_1
        \big)}
\newline
\hspace*{0.cm}    -H_{\lambda_2} \big(
                \sum_{i=1}^{\lambda_2} \frac{\binom{i
                +\lambda_1
                }{i}}{2^{i}}-1\big)
        +
        \sum_{i=1}^{\lambda_2}\tfrac{\tbinom{i
        +\lambda_1
        }{i}H_{i
        +\lambda_1
        }}{2^{i}}
\Bigg)
+\frac{\binom{\lambda_1
+\lambda_2
}{\lambda_2} \big(
        -H_{\lambda_1
        -\lambda_2
        }
        +H_{\lambda_1}
\big)}{8(\lambda_1+1)}
\big(3 \lambda_1
        +4 \lambda_1^2
        +\lambda_1^3
        -3 \lambda_2
        +5 \lambda_1 \lambda_2
        -9 \lambda_2^2
        -5 \lambda_1 \lambda_2^2
        +4 \lambda_2^3
\big)-2^{\lambda_2-5} \big(
        14
        -3 \lambda_1
        +3 \lambda_1^2
        +\lambda_2
        -6 \lambda_1 \lambda_2
        +3 \lambda_2^2
\big)
\big(
        \sum_{i=1}^{\lambda_2} 2^{-i} \binom{i
        +\lambda_1
        }{i}+1\big)
\newline+\frac{\binom{\lambda_1
+\lambda_2
}{\lambda_2}}{\big(
        \lambda_1+1\big)^2} \frac{1}{32} \big(
        14
        +3 \lambda_1^4
        +\lambda_1^3 \big(
                3-16 \lambda_2\big)
        -5 \lambda_2
        -51 \lambda_2^2
        +24 \lambda_2^3
        +\lambda_1^2 \big(
                5 \lambda_2^2-5 \lambda_2+11\big)
        +\lambda_1 \big(
                8 \lambda_2^3-30 \lambda_2^2+38 \lambda_2+25\big)
\big)
\newline +\frac{\binom{\lambda_1
+\lambda_2
}{\lambda_2}}{\lambda_1+1}
 \frac{1}{16} \big(
        8
        +5 \lambda_1
        -4 \lambda_1^2
        -\lambda_1^3
        +7 \lambda_2
        -9 \lambda_1 \lambda_2
        +17 \lambda_2^2
        +9 \lambda_1 \lambda_2^2
        -8 \lambda_2^3
\big) H_{\lambda_1
+\lambda_2
}\\
\end{mma}

\smallskip

We emphasize that each calculation step of this transformation can be verified in the same fashion as it has been worked out for the sum $S_0(\lambda_1,\lambda_2)$ from above.
Completely, analogously we obtain the representations \texttt{solSum2}, \texttt{solSum3}, \texttt{solSum4}, \texttt{solSum5}  (together with correctness proofs) of the remaining four sums in~\eqref{eq:sixdouble} in terms of indefinite nested sums over hypergeometric products which are all valid for all $\lambda_1,\lambda_2$ with $1\leq\lambda_2\leq\lambda_1$. This completes step \textbf{(DEF)} of our summation tactic.

Now we are ready to carry out step~\textbf{(IND)}.
With the function call
\begin{mma}
\In \{solSum0,solSum1,solSum2,solSum3,solSum4,solSum5\}=\newline
\hspace*{2cm}SigmaReduce[\{solSum0,solSum1,solSum2,solSum3,solSum4,solSum5\},\lambda_2];\\ 
\end{mma}
\noindent we synchronize the arising sums and products within the representations of the six sums. More precisely, after this transformation all expressions depend only on the following indefinite nested sums defined over hypergeometric products
\begin{multline}\label{Equ:Equ:ProdSumRep}
2^{\lambda_2},\binom{\lambda_1
+\lambda_2
}{\lambda_2},H_{\lambda_1},H_{\lambda_1
-\lambda_2
},H_{\lambda_2},H_{\lambda_1
+\lambda_2
},
\sum_{i=1}^{\lambda_2} 2^{-i} \binom{i
+\lambda_1
}{i},
\sum_{i=1}^{\lambda_2} 2^{-i} \binom{i
+\lambda_1
}{i} H_{-i
+\lambda_1
},\\
\sum_{i=1}^{\lambda_2} 2^{-i} \binom{i
+\lambda_1
}{i} H_{i
+\lambda_1
},
\sum_{i=1}^{\lambda_2} \frac{2^i}{i \binom{i
+\lambda_1
}{i}},
\sum_{i=1}^{\lambda_2} \frac{2^i 
\sum_{j=1}^i 2^{-j} \binom{j
+\lambda_1
}{j}}{i \binom{i
+\lambda_1
}{i}},
\sum_{i=1}^{\lambda_2} \frac{2^i \big(
        \sum_{j=1}^i 2^{-j} \binom{j
        +\lambda_1
        }{j}\big)^2}{\binom{i
+\lambda_1
}{i} \big(
        1
        +i
        +\lambda_1
\big)}
\end{multline}
whose sequences produced by the objects in~\eqref{Equ:Equ:ProdSumRep} are algebraically independent among each other.

Putting all the building blocks together 
\begin{mma}
 \In \tilde{T}=\frac{\lambda_1(\lambda_1-1)}{4}+\frac{\lambda_2(\lambda_2-3)}{4}-2\,solSum0-
 \Bigg(\left(\binom{\lambda_1}{2}+\binom{\lambda_2+1}{2}\right)(H_{\lambda_1+\lambda_2}-1)+\newline\frac{\lambda_1-\lambda_2+1}{\lambda_1+1}\binom{\lambda_1+\lambda_2}{\lambda_2}\Big(
 -solSum1+solSum2-solSum3-solSum4+solSum5\Big)\Bigg);\\
\end{mma}
\noindent yields an alternative representation $\tilde{T}$ of our expression $T(\lambda_1,\lambda_2)$ given in terms of~\eqref{Equ:ExprAMinusB} which is valid for all $\lambda_1,\lambda_2\in\N$ with $1\leq\lambda_2\leq\lambda_1$. In a nutshell, we solved Problem~S and obtained the simplification $\tilde{T}$ where all arising sums and products are algebraically independent among each other. Finally, we
consider all the sums and products in this expression as variables, put them over a common denominator and expand the derived numerator. More precisely, we apply the standard Mathematica command \texttt{Together} to $\tilde{\texttt{T}}$ and obtain by simple polynomial arithmetic the answer 
\begin{mma}
\In Together[\tilde{T}]\\
\Out 0\\
\end{mma}
\noindent This implies that $T(\lambda_1,\lambda_2)=0$ holds or equivalently that~\eqref{Equ:A=B} holds 
for all $\lambda_1,\lambda_2\in\N$ with $1\leq\lambda_2\leq\lambda_1$. This completes the proof of \reft{zweizeiler}.

As a reward for all our calculations we obtain besides a proof of Conjecture~\ref{Conjecture:zweizeiler}
in addition new representations
of $C(\lambda)$.  Using~\eqref{eq:zweizeiler} together with~\eqref{Equ:S0Rep} yields
\begin{multline*}
C(\lambda)=C(\lambda_1,\lambda_2)=\frac{\lambda_1(\lambda_1-1)}{4}+\frac{\lambda_2(\lambda_2-3)}{4}\\[-0.5cm]
-2\Big(-\frac{2^{\lambda_2}}{\binom{\lambda_1
+\lambda_2
}{\lambda_2}} \Big(
        \sum_{i=1}^{\lambda_2} \frac{\binom{i
        +\lambda_1
        }{i}}{2^{i}}+1\Big)
+\frac{-1
        -\lambda_1
        +\lambda_2}{2}  
\sum_{i=1}^{\lambda_2} \frac{\displaystyle 2^i 
\sum_{j=1}^i\tfrac{\binom{j
+\lambda_1
}{j}}{2^{j}}}{i \binom{i
+\lambda_1
}{i}}\\
-\frac{1
        +\lambda_1
        -\lambda_2}{2}
\sum_{i=1}^{\lambda_2} \frac{2^i}{i \binom{i
+\lambda_1
}{i}}
+\big(  1
        +\lambda_1
        -\lambda_2
\big)\big( 
H_{\lambda_1}
+\frac{1}{2}H_{\lambda_2}-H_{\lambda_1
-\lambda_2
}\big)
+1\Big).
\end{multline*}
This expression is of particular interest if one wants to calculate $C(\lambda_1,\lambda_2)$ efficiently for $\lambda_2=0,1,2,3,\dots$ and keeping $\lambda_1$ symbolic. 
In addition, we obtain the specially nice formula
$$C(\lambda_2,\lambda_2)=
\frac{\lambda_2(\lambda_2-2)}{2}
-2\Big(-\frac{2^{2 \lambda_2}}{\binom{2 \lambda_2}{\lambda_2}}
+\frac{H_{\lambda_2}}{2}
+1\Big)$$
using the identity~\eqref{Equ:S0SpecialCase}. Similarly, if one is interested in an efficient evaluation of $C(\lambda)$ for a symbolic $\lambda_2$ and a fixed distance $\delta=\lambda_1-\lambda_2\geq0$, we can derive with our summation toolbox the following representation for the sum $S_0(\lambda_1,\lambda_2)$:
\begin{mma}
\In resDelta=EvaluateMultiSum[\frac{(-1)^k \binom{\lambda_2}{k} (-2+2 k)!}{\big(
        2
        +\lambda_1
        -\lambda_2
\big)_{-1+2 k}}/.\lambda_1\to\delta+\lambda_2,\newline
\hspace*{7cm}\{\{k,1,\lambda_2\}\},\{\delta,\lambda_2\},\{0,1\},\{\infty,\infty\}];\\ 
\end{mma} 
 
\begin{mma}
\In SigmaReduce[resDelta,\delta,Tower\to \big\{\binom{\delta 
+\lambda_2
}{\delta },\binom{\delta 
+2 \lambda_2
}{\delta },H_{\delta 
+\lambda_2
},H_{\delta 
+2 \lambda_2
},H_{\delta }\big\}] \\
\Out \big(
        -1
        -\delta 
        +
        \frac{(\delta +1)2^{2 \lambda_2}}{\binom{2 \lambda_2}{\lambda_2}}
\big) 
\sum_{i=1}^{\delta } \frac{2^i \binom{i
+\lambda_2
}{i}}{\binom{i
+2 \lambda_2
}{i} \big(
        1
        +i
        +2 \lambda_2
\big)}
+\frac{2^{\delta +2} \lambda_2\big(
        1
        +\delta 
        +\lambda_2
\big) \binom{\delta 
+\lambda_2
}{\delta }}{\binom{\delta 
+2 \lambda_2
}{\delta } \big(
        1
        +\delta 
        +2 \lambda_2
\big)} \sum_{i=1}^{\delta } \frac{2^{-i} \binom{i
        +2 \lambda_2
        }{i}}{\binom{i
        +\lambda_2
        }{i} \big(
                i
                +2 \lambda_2
        \big)}
\newline
-2 (\delta +1)\lambda_2
        \sum_{i=1}^{\delta } \frac{2^i \binom{i
        +\lambda_2
        }{i}\displaystyle 
        \sum_{j=1}^i \frac{2^{-j} \binom{j
        +2 \lambda_2
        }{j}}{\binom{j
        +\lambda_2
        }{j} \big(
                j
                +2 \lambda_2
        \big)}}{\binom{i
        +2 \lambda_2
        }{i} \big(
                1
                +i
                +2 \lambda_2
        \big)}
+\frac{2^{2 \lambda_2}(\delta +1)}{\binom{2 \lambda_2}{\lambda_2} \big(
        2 \lambda_2+1\big)}
+\frac{2^{\delta +1} \big(1
        +\delta 
        +\lambda_2
\big)\binom{\delta 
+\lambda_2
}{\delta }}{\binom{\delta 
+2 \lambda_2
}{\delta } \binom{2 \lambda_2}{\lambda_2}\big(
        1
        +\delta 
        +2 \lambda_2
\big)}\times\newline 
\times \big(
        -2^{2 \lambda_2}
        +\binom{2 \lambda_2}{\lambda_2}
\big)
-\frac{\delta +1}{2 \lambda_2+1}
+(\delta +1) H_{\delta 
+2 \lambda_2
}
+\frac{1}{2} (\delta +1) H_{\lambda_2}
-(\delta +1) H_{2 \lambda_2}
-(\delta +1) H_{\delta }\\
\end{mma}
\noindent This yields
\begin{align*}
C(\lambda)=C(\lambda_1,\lambda_2)&=\frac{\lambda_1(\lambda_1-1)}{4}+\frac{\lambda_2(\lambda_2-3)}{4}\\
-2\Bigg(&\big(
        -1
        -\delta 
        +
        \frac{(\delta +1)2^{2 \lambda_2}}{\binom{2 \lambda_2}{\lambda_2}}
\big) 
\sum_{i=1}^{\delta } \frac{2^i \binom{i
+\lambda_2
}{i}}{\binom{i
+2 \lambda_2
}{i} \big(
        1
        +i
        +2 \lambda_2
\big)}\\
&+\frac{2^{\delta +2} \lambda_2\big(
        1
        +\delta 
        +\lambda_2
\big) \binom{\delta 
+\lambda_2
}{\delta }}{\binom{\delta 
+2 \lambda_2
}{\delta } \big(
        1
        +\delta 
        +2 \lambda_2
\big)} \sum_{i=1}^{\delta } \frac{2^{-i} \binom{i
        +2 \lambda_2
        }{i}}{\binom{i
        +\lambda_2
        }{i} \big(
                i
                +2 \lambda_2
        \big)}\\
&-2 (\delta +1)\lambda_2
        \sum_{i=1}^{\delta } \frac{2^i \binom{i
        +\lambda_2
        }{i}\displaystyle 
        \sum_{j=1}^i \frac{2^{-j} \binom{j
        +2 \lambda_2
        }{j}}{\binom{j
        +\lambda_2
        }{j} \big(
                j
                +2 \lambda_2
        \big)}}{\binom{i
        +2 \lambda_2
        }{i} \big(
                1
                +i
                +2 \lambda_2
        \big)}\\
&+\frac{2^{2 \lambda_2}(\delta +1)}{\binom{2 \lambda_2}{\lambda_2} \big(
        2 \lambda_2+1\big)}
+\frac{2^{\delta +1} \big(1
        +\delta 
        +\lambda_2
\big)\binom{\delta 
+\lambda_2
}{\delta }}{\binom{\delta 
+2 \lambda_2
}{\delta } \binom{2 \lambda_2}{\lambda_2} \big(
        1
        +\delta 
        +2 \lambda_2
\big)} \big(
        -2^{2 \lambda_2}
        +\binom{2 \lambda_2}{\lambda_2}
\big)\\
&-\frac{\delta +1}{2 \lambda_2+1}
+(\delta +1) H_{\delta 
+2 \lambda_2
}
+\frac{1}{2} (\delta +1) H_{\lambda_2}
-(\delta +1) H_{2 \lambda_2}
-(\delta +1) H_{\delta }\Bigg)
\end{align*}
with $\lambda_1=\lambda_2+\delta$ for a positive integer parameter $\lambda_2$ and a non-negative integer $\delta$.

\subsection*{Acknowledgements}

The authors thank Christian Krattenthaler for introducing them to the problems treated in this paper, and the anonymous referee for the useful suggestions.



\newcommand{\etalchar}[1]{$^{#1}$}

\end{document}